\documentstyle[12pt]{article}

\newcounter{primer}
\newcounter{tebr}
\newcounter{lebr}

\def\teo{\addtocounter{tebr}{1}{\vspace{.2cm}\noindent T{\footnotesize HEOREM}
\arabic{section}.\arabic{tebr}.}\hspace{1em}}
\def\lema{\addtocounter{lebr}{1}{\vspace{.2cm}\noindent L{\footnotesize EMMA}
\arabic{section}.\arabic{lebr}.}\hspace{1em}}
\def\exmp{\addtocounter{primer}{1}{\vspace{.2cm}\noindent E{\footnotesize
XAMPLE} \arabic{section}.\arabic{primer}.}\hspace{1em}}
\newcommand{\qed}{\hfill $\Box$ \\[0.2cm]}
\newcommand{\dkz}{\vspace{.12cm}\noindent{P{\footnotesize ROOF}.}\hspace{1em} \nopagebreak}

\def\por#1{\makebox[3em][l]{$#1$}}

\def\ppor#1{\makebox[20em][l]{$#1$}}

\def\pppor#1{\makebox[40em][l]{$#1$}}

\def\mix#1{\begin{picture}(10,10)\put(5,3){\circle{10}}
\put(5,3){\makebox(0,0){\tiny$#1$}}\end{picture}}

\def\d#1{\mbox{\tiny +#1}}

\def\l#1{\mbox{\tiny #1+}}

\newcommand{\gdna}{
\begin{picture}(24,12)(0,6)
\put(0,10){\vector(1,0){23}}
\put(9,14){\circle*{2}}
\put(14,14){\circle*{2}}
\put(12,5){\makebox(0,0){\scriptsize$\Gamma$}}
\end{picture}}

\newcommand{\fdna}{
\begin{picture}(24,12)(0,6)
\put(0,10){\vector(1,0){23}}
\put(9,14){\circle*{2}}
\put(14,14){\circle*{2}}
\put(12,5){\makebox(0,0){\scriptsize$\Phi$}}
\end{picture}}

\newcommand{\pdna}{
\begin{picture}(24,12)(0,6)
\put(0,10){\vector(1,0){23}}
\put(9,14){\circle*{2}}
\put(14,14){\circle*{2}}
\put(12,5){\makebox(0,0){\scriptsize$\Psi$}}
\end{picture}}

\newcommand{\pj}{1\!\!1}
\newcommand{\uf}{[\,\mf]}

\newcommand{\md}{\mbox{\boldmath{${\delta}$}}}

\newcommand{\mep}{\mbox{\boldmath{${\varepsilon}$}}}
\newcommand{\met}{\mbox{\boldmath{${\eta}$}}}
\newcommand{\mfv}{\mbox{\boldmath{$F$}}}
\newcommand{\mgv}{\mbox{\boldmath{$G$}}}
\newcommand{\mtv}{\mbox{\boldmath{$T$}}}
\newcommand{\msv}{\mbox{\boldmath{$S$}}}
\newcommand{\mrv}{\mbox{\boldmath{$R$}}}
\newcommand{\mj}{\mbox{\bf 1}}
\newcommand{\ccb}{$Constr({\cal B})$}
\newcommand{\cb}{{\cal B}}
\newcommand{\bcc}{\mbox{\bf BiCartCl}}
\newcommand{\mb}{\mbox{\bf b}}
\newcommand{\mc}{\mbox{\bf c}}
\newcommand{\mw}{\mbox{\bf w}}
\newcommand{\mk}{\mbox{\bf k}}
\newcommand{\ml}{\mbox{\boldmath{$l$}}}
\newcommand{\mm}{\mbox{\bf m}}
\newcommand{\mf}{\mbox{\bf f}\,}
\newcommand{\ri}{\mbox{\rm I}}
\newcommand{\mi}{\mbox{\scriptsize{\rm I}}}
\newcommand{\ro}{\mbox{\rm O}}

\newcommand{\str}{\rightarrow}
\newcommand{\rts}{\leftarrow}

\newcommand{\HDS}{\vrule width0pt height2.3ex depth1.05ex\displaystyle}
\newcommand{\J}{$\cal J$}

\newcommand{\F}{${\cal F}$}
\newcommand{\I}{\wedge}
\newcommand{\IL}{\vee}

\def\f#1#2{{{\HDS #1}\over{\HDS #2}}}

\def\fp#1#2{{{\HDS #1}\atop{\HDS #2}}}

\def\afrac#1{{\phantom{\HDS #1}\atop{\HDS #1}}}

\def\bfrac#1#2{{\phantom{\HDS \Theta}\atop{\phantom{
\vrule width0pt height2.3ex depth#2\displaystyle \Theta}\atop{\HDS #1}}}}

\def\rza{{\mbox{\hspace{0.5em}}}}

\def\rzb{{\mbox{\hspace{2em}}}}

\def\rzc{{\mbox{\hspace{7em}}}}

\def\Bpak#1#2#3{\f{\f{#1}{\cdots} \pravila{#2}}{#3}}

\def\pravilo#1{ \rza \makebox[-.5em][l]{\mbox{\rm #1}}}

\def\pravila#1{ \rza \makebox[-.5em][l]{\raisebox{-1.7ex}
{\parbox{15ex}{\footnotesize \baselineskip=0.5\baselineskip #1}}}}

\title{G-dinaturality}

\author{\sc Zoran Petri\' c\\
Matemati\v cki institut\\
SANU, Knez Mihailova 35, p.f. 367\\
email: zpetric@mi.sanu.ac.yu}

\date{\today}

\begin{document}

\maketitle

\abstract{An extension of the notion of dinatural transformation is introduced
in order to give a criterion for preservation of dinaturality
under composition. An example of an application is given by
proving that all bicartesian closed canonical transformations are dinatural.
An alternative sequent system for
intuitionistic propositional logic is introduced as a device, and a cut
elimination procedure is established for this system.}

\section{Introduction}

The aim of tis paper is to introduce a generalization of the traditional
notion of dinaturality and to give a geometrical criterion for preservation
of dinaturality under composition.
A certain importance is usually attached to this question (see \cite{EK66},
\cite{BFSS90}, \cite{GSS92} and \cite{B93}). It is useful to consult
\cite{M76}
to find about the historical perspective of the notion of naturality.
The first
extension of this notion towards our g-dinatural
transformations was given by Eilenberg and Kelly in \cite{EK66}. In the
present paper
we generalize the definition of dinaturality introduced by
Dubuc and Street in
\cite{DS70}.
The theory of g-dinatural transformations is here applied to 
bicartesian closed canonical transformations. It is proved that they
are all dinatural in the sense of \cite{DS70}.

By a bicartesian closed category we mean a category equipped with finite products
and coproducts, including initial and terminal objects, which is closed
in the sense that for every object $A$, the functor $A\times\_$ has the right
adjoint $A\str\_$. This category may serve as a framework for the categorial
proof theory of intuitionistic propositional logic. However, despite
that we are using a very traditional categorial object, our notation and
definitions are a little bit unusual. This choice is forced by the technique
that we intend to use here, and we believe that it is optimal.
\\[.2cm]
{\bf Notation.} For objects we use the schematic letters $A,B,C,\ldots,
A_1,\ldots$, and for morphisms the schematic letters $f,g,h,\ldots,f_1,
\ldots$ The product of $A$ and $B$ is denoted by $A\times B$, and the
coproduct by $A+B$. We use O and I to specify the initial and the terminal
object of a category. To denote that a morphism $f$ has the source $A$
and the
target $B$ we use the notation $f:A\vdash B$, and we say that $f$ is of the
type $A\vdash B$. Apart from the logical motivation for the symbol $\vdash$
instead of $\str$, we have another reason, which comes from our intention
to write complex objects linearly: we use $A\str B$ instead of
$B^{\mbox{$A$}}$ for exponentiation, i.e. the immage of $B$ under the
right adjoint of the functor $A\times\_$. However, in Section 2, where we
deal with a new notion of dinaturality, and which is, except for examples,
self-contained, we use the standard symbol $\str$ for morphisms. Also,
to avoid too many parentheses, we assume that the morphism operation
$\circ$ binds more strongly than $\times,+,\str$; for example, we write
$g\circ f\times h$ for $(g\circ f)\times h$.
\\[.2cm]
{\bf Bicartesian closed categories.} A bicartesian closed category
$\cal B$ satisfies:
\\[.1cm]
For every triple $A,B,C$ of objects from $\cal B$, we have the following
{\em special morphisms} in $\cal B$
\[\begin{array}{ll}
\mj_A:A\vdash A, &
\\[.1cm]
\md_A:A\times I\vdash A, & \md^i_A:A\vdash A\times I,
\\[.1cm]
\ppor{\mb^\str_{A,B,C}:A\times(B\times C)\vdash (A\times B)\times C,}
& \ppor{\mb^\rts_{A,B,C}:(A\times B)\times C\vdash A\times (B\times C),}
\\[.1cm]
\mc_{A,B}:A\times B\vdash B\times A, &
\\[.1cm]
\mw_A:A\vdash A\times A, &
\mk_A:A\vdash \ri, 
\\[.1cm]
\mm_A:A+A\vdash A,  & \ml_A:\ro\vdash A,
\\[.1cm]
\ml^1_{A,B}:A\vdash A+B, & \ml^2_{A,B}:B\vdash A+B,
\\[.1cm]
\mep_{A,B}:A\times(A\str B)\vdash B,
& \met_{A,B}:B\vdash A\str(A\times B),
\end{array}\]
and the following operations on morphisms:
\[\begin{array}{ll}
\ppor{\f{f:A\vdash B \rzb g:B\vdash C}{g\circ f:A\vdash C},} &
\ppor{\f{f:A\vdash B \rzb g:C\vdash D}{f\str g:B\str C\vdash A\str D},}
\\[.3cm]
\f{f:A\vdash B \rzb g:C\vdash D}{f\times g:A\times C\vdash B\times D},\rzb\rzb &
\f{f:A\vdash B \rzb g:C\vdash D}{f+ g:A+ C\vdash B+ D}.
\end{array}\]
Also, the following equations must be satisfied
\[\begin{array}{ll}
\por{(cat\; 1)} & \pppor{\mj_B\circ f=f\circ\mj_A=f,}
\\[.1cm]
(cat\; 2) & h\circ(g\circ f)=(h\circ g)\circ f,
\end{array}\]
\[\begin{array}{ll}
\por{(\times 1)} & \pppor{\mj_A\times \mj_B=\mj_{A\times B},}
\\[.1cm]
(\times 2) &  (g_1\circ g_2)\times (f_1\circ f_2)=(g_1\times f_1)\circ
(g_2\times f_2),
\\[.1cm]
(\md) &  f\circ\md_A=\md_B\circ(f\times\mj_{\mi})
\\[.1cm]
(\md\md^i) &  \md_A\circ\md^i_A=\mj_A,\quad\quad\quad
\md^i_A\circ\md_A=\mj_{A\times{\mi}},
\\[.1cm]
(\md\mc) &  \md_{\mi}\circ\mc=\md_{\mi},
\\[.1cm]
(\mb) &  ((f\times g)\times h)\circ\mb^\str_{A,B,C}=
\mb^\str_{D,E,F}\circ(f\times(g\times h)),
\\[.1cm]
(\mb\mb) &  \mb^\str_{A,B,C}\circ\mb^\rts_{A,B,C}=\mj_{(A\times B)\times C},
\quad\quad\quad
\mb^\rts_{A,B,C}\circ\mb^\str_{A,B,C}=\mj_{A\times (B\times C)},
\\[.1cm]
(\mb 5) &  \mb^\str_{A\times B,C,D}\circ\mb^\str_{A,B,C\times D}=
(\mb^\str_{A,B,C}\times\mj_D)\circ\mb^\str_{A,B\times C,D}\circ
(\mj_A\times\mb^\str_{B,C,D}),
\\[.1cm]
(\mc) &  (g\times f)\circ\mc_{A,B}=\mc_{C,D}\circ(f\times g),
\\[.1cm]
(\mc\mc) &  \mc_{B,A}\circ\mc_{A,B}=\mj_{A\times B},
\\[.1cm]
(\mb\mc\md) &  (\md_A\times\mj_B)\circ\mb^\str_{A,{\mi},B}=
\mj_A\times\md_B\circ\mc_{{\mi},B},
\\[.1cm]
(\mb\mc 6) &  \mb^\str_{C,A,B}\circ\mc_{A\times B,C}\circ\mb^\str_{A,B,C}=
(\mc_{A,C}\times\mj_B)\circ\mb^\str_{A,C,B}\circ(\mj_A\times\mc_{B,C}),
\end{array}\]
\[\begin{array}{ll}
\por{(\mw)} & \pppor{(f\times f)\circ\mw_A=\mw_B\circ f,}
\\[.1cm]
(\md\mw) &  \md_{\mi}\circ\mw_{\mi}=\mj_{\mi},
\\[.1cm]
(\mb\mw) &  \mb^\str_{A,A,A}\circ(\mj_A\times\mw_A)\circ\mw_A=
(\mw_A\times\mj_A)\circ\mw_A,
\\[.1cm]
(\mc\mw) &  \mc_{A,A}\circ\mw_A=\mw_A,
\\[.1cm]
(\mb\mc\mw 8) &  \mc^m_{A,B,A,B}\circ\mw_{A\times B}=\mw_A\times\mw_B,\quad
{\mbox{\rm where}}
\\
\por{\hspace{2em}\mc^m_{A,B,C,D}=^{df} \mb^\str_{A,C,B\times D}
\circ(\mj_A\times(\mb^\rts_{C,B,D}\circ(\mc_{B,C}\times\mj_D)\circ
\mb^\str_{B,C,D}))\circ\mb^\rts_{A,B,C\times D},} &
\\[.1cm]
(\mk) &  {\mbox{\rm for }}\; f:A\vdash \ri,\; f=\mk_A,
\\[.1cm]
(\md\mk\mw) &  \md_A\circ(\mj_A\times\mk_A)\circ\mw_A=\mj_A,
\end{array}\]
\[\begin{array}{ll}
\por{(+ 1)} &  \pppor{\mj_A +\mj_B=\mj_{A+B},}
\\[.1cm]
(+ 2) &  (g_1\circ g_2)+(f_1\circ f_2)=(g_1+ f_1)\circ(g_2 + f_2),
\\[.1cm]
(\ml^1) &  (f_1+f_2)\circ\ml^1_{A_1,A_2}=\ml^1_{B_1,B_2}\circ f_1,
\\[.1cm]
(\ml^2) &  (f_1+f_2)\circ\ml^2_{A_1,A_2}=\ml^2_{B_1,B_2}\circ f_2,
\\[.1cm]
(\ml) &  {\mbox{\rm for }}\; f:\ro\vdash A,\; f=\ml_A,
\\[.1cm]
(\mm) &  f\circ\mm_A=\mm_B\circ(f+f),
\\[.1cm]
(\ml\mm 1) &  \mm_A\circ\ml^1_{A,A}=\mj_A=\mm_A\circ\ml^2_{A,A},
\\[.1cm]
(\ml\mm 2) &  \mm_{A+B}\circ(\ml^1_{A,B}+\ml^2_{A,B})=\mj_{A+B},
\end{array}\]
\[\begin{array}{ll}
\por{(\str 1)} &  \pppor{\mj_A\str\mj_B=\mj_{A\str B},}
\\[.1cm]
(\str 2) &  (g_1\circ g_2)\str(f_1\circ f_2)=(g_2\str f_1)\circ(g_1\str f_1),
\\[.1cm]
(\mep 1) &  f\circ\mep_{C,A}=\mep_{C,B}\circ(\mj_C\times(\mj_C\str f)),
\\[.1cm]
(\met 1) &  (\mj_C\str(\mj_C\times f))\circ \met_{C,A}=\met_{C,B}\circ f,
\\[.1cm]
(\mep 2) &  \mep_{B,C}\circ(f\times(\mj_B\str\mj_C))=\mep_{A,C}\circ
(\mj_A\times(f\str\mj_C)),
\\[.1cm]
(\met 2) &  (f\str(\mj_B\times\mj_C))\circ\met_{B,C}=
(\mj_A\str(f\times\mj_C))\circ\met_{A,C},
\\[.1cm]
(\mj\mep\met) &  (\mj_A\str\mep_{A,B})\circ\met_{A,A\str B}=\mj_{A\str B},
\\[.1cm]
(\mep\mj\met) &  \mep_{A,A\times B}\circ(\mj_A\times\met_{A,B})=\mj_{A\times B}.
\end{array}\]
It is easy to extract the definitions of {\em symmetric monoidal closed},
{\em cartesian closed} and {\em bicartesian} categories from the definition
above. This is the first reason to accept the approach above to bicartesian
closed categories. Another reason is the sequent system that we are going to
deal with, and the process of cut elimination tied to it.

The proof that the above definition is equivalent to the equational definition
of bicartesian closed categories given in \cite{LS86} requires some effort,
but we won't go into this matter here.

\section{Graphs and g-dinatural transformations}

This section is devoted to the notion of g-dinatural transformations.
These transformations will serve as morphisms in a functor category whose
object are functors of the type ${\cal A}_1\times{\cal A}_2\times\ldots
{\cal A}_n\str{\cal A}$ for some category $\cal A$, ${\cal A}_i\in\{{\cal A},
{\cal A}^{op}\}$ and arbitrary $n\in N$. They are always equipped with
``graphs'' and this explains the letter g in the label g-dinatural.
First we define the notion of graph. 

For  $m\geq 0$ and $n\geq 0$, let
$M$ be the set $\{x_1,x_2,\ldots,x_m\}+\{y_1,y_2,\ldots , y_n\}$
whose elements
we call {\em argument places} (the $x$'s are {\em left-hand side}
argument places and the
$y$'s are {\em right-hand side} argument places, and if $m=n=0$,
$M$ is empty). Let $G$ be a finite set and let $l$
be a mapping $l:M \str \{1,-1\}$, which intuitively denotes the covariance or
the contravariance of an argument place. If $l(u)=1$ we say that $u$ is
a {\em positive argument place} and we write $u^+$,
and if $l(u)=-1$ we call it {\em negative argument place} and we write $u^-$.
The elements of the set $V=M \cup G$ are called {\em vertices}.

Let $E$ be a set of pairs of elements from $V$ that we call {\em edges}.
Let $u \sim v$ means that there is an edge $\{u,v\}$, and let $\simeq$
be the reflexive and transitive closure of $\sim$. Then the
equivalence classes of $\simeq$ together with the corresponding
edges from $E$, are called {\em components}. Let us
enumerate these components by $1,2,\ldots , k$, ($k\geq 0$), and let
$\pi$ be the mapping $\pi: V\str \{1,\ldots ,k\}$ that maps a vertex from $V$
to the number of its component. We call this function
{\em component classifier}.

For $V$, $l$ and $E$ as above, the triple $(V,l,E)$ is called {\em graph}
iff the following conditions hold:
\begin{enumerate}
\item every vertex belongs to some edge,
\item $\{x_i,x_j\}\in E$\quad iff\quad $l(x_i)=-l(x_j)$
and $x_i,x_j$ are in the same component,
\item $\{y_i,y_j\}\in E$\quad iff\quad $l(y_i)=-l(y_j)$
and $y_i,y_j$ are in the same component,
\item $\{x_i,y_j\}\in E$\quad iff\quad $l(x_i)=l(y_j)$
and $x_i,y_j$ are in the same component,
\item if a component $K$ includes an edge between two
argument places, then $K \cap G =\emptyset$; otherwise, $K\cap G=\{g\}$
for some vertex $g\in G$, and for every $u\in K\setminus\{g\}$ the edge
$\{u,g\}$ is in $E$ ($K$).
\end{enumerate}

\exmp The following diagram illustrates a graph with 3 components,
wher $G$ is a singleton.

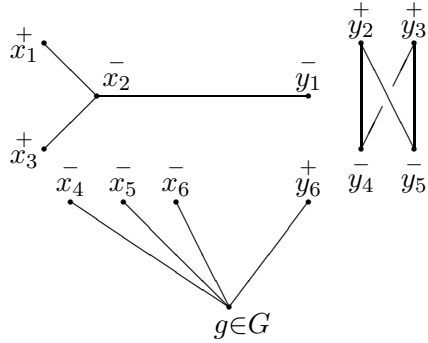
\begin{center}
\begin{picture}(180,140)
\put(20,120){\circle*{2}}
\put(18,120){\makebox(0,0)[r]{$\stackrel{+}{x_{1}}$}}
\put(20,80){\circle*{2}}
\put(18,80){\makebox(0,0)[r]{$\stackrel{+}{x_{3}}$}}
\put(40,100){\circle*{2}}
\put(42,104){\makebox(0,0)[bl]{$\stackrel{-}{x_{2}}$}}
\put(30,60){\circle*{2}}
\put(30,64){\makebox(0,0)[b]{$\stackrel{-}{x_{4}}$}}
\put(50,60){\circle*{2}}
\put(50,64){\makebox(0,0)[b]{$\stackrel{-}{x_{5}}$}}
\put(70,60){\circle*{2}}
\put(70,64){\makebox(0,0)[b]{$\stackrel{-}{x_{6}}$}}
\put(120,60){\circle*{2}}
\put(120,64){\makebox(0,0)[b]{$\stackrel{+}{y_{6}}$}}
\put(120,100){\circle*{2}}
\put(120,104){\makebox(0,0)[b]{$\stackrel{-}{y_{1}}$}}
\put(140,120){\circle*{2}}
\put(140,124){\makebox(0,0)[b]{$\stackrel{+}{y_{2}}$}}
\put(160,120){\circle*{2}}
\put(160,124){\makebox(0,0)[b]{$\stackrel{+}{y_{3}}$}}
\put(140,80){\circle*{2}}
\put(140,78){\makebox(0,0)[t]{$\stackrel{-}{y_{4}}$}}
\put(160,80){\circle*{2}}
\put(160,78){\makebox(0,0)[t]{$\stackrel{-}{y_{5}}$}}
\put(90,20){\circle*{2}}
\put(95,18){\makebox(0,0)[t]{$g{\in}G$}}
\put(20,120){\line(1,-1){20}}
\put(20,80){\line(1,1){20}}
\put(40,100){\line(1,0){80}}
\put(140,80){\line(0,1){40}}
\put(140,80){\line(1,2){8}}
\put(160,120){\line(-1,-2){8}}
\put(160,80){\line(0,1){40}}
\put(160,80){\line(-1,2){20}}
\put(90,20){\line(-3,2){60}}
\put(90,20){\line(-1,2){20}}
\put(90,20){\line(-1,1){40}}
\put(90,20){\line(3,4){30}}
\end{picture}
\end{center}

\vspace{.2cm}
For a graph $\Gamma$ we denote by $\Gamma_i$ its $i$-th connectional
component. Let $\Gamma_i^+$ and $\Gamma_i^-$ be the sets of positive and
negative vertices from $\Gamma_i$ respectively.
Also, for a symbol $a$ we use the abbreviation
$a^k$ for the sequence of $k$ occurrences of this symbol.

Let for a single component graph $\Gamma=(V,l,E)$ and a small category
$\cal A$, $T$ and $S$ be functors
\[\begin{array}{l}
T:{\cal A}^{l(x_1)} \times\ldots\times {\cal A}^{l(x_m)} \str {\cal A},\\[.2cm]
S:{\cal A}^{l(y_1)} \times\ldots\times {\cal A}^{l(y_n)} \str {\cal A},
\end{array}\]
where ${\cal A}^1 ={\cal A}$ and ${\cal A}^{-1} ={\cal A}^{op}$.
Let $\alpha$ be a set
\[\alpha=\{\alpha(A):T(A^m)\str S(A^n) \; | \; A\in Ob({\cal A})\}\]
of morphisms from $\cal A$ indexed by the set of objects from $\cal A$.
Such a family is called a {\em transformation}. Then we say that $\alpha$ is
a {\em g-dinatural transformation} from $T$ to $S$ with the graph $\Gamma$,
which is denoted by $\alpha: T \gdna S$, if for every pair of
objects $A$, $C$ and
every morphism $f:A \str C$ from $\cal A$, the following diagram commutes:

\begin{center}
\begin{picture}(200,110)
\put(80,20){\vector(1,0){40}}
\put(80,90){\vector(1,0){40}}
\put(30,50){\vector(3,-2){30}}
\put(140,80){\vector(3,-2){30}}
\put(30,60){\vector(3,2){30}}
\put(140,30){\vector(3,2){30}}
\put(78,20){\makebox(0,0)[r]{$T\langle A,A\rangle$}}
\put(78,90){\makebox(0,0)[r]{$T\langle C,C\rangle$}}
\put(122,20){\makebox(0,0)[l]{$S\langle A,A\rangle$}}
\put(122,90){\makebox(0,0)[l]{$S\langle C,C\rangle$}}
\put(25,55){\makebox(0,0)[r]{$T\langle A,C\rangle$}}
\put(175,55){\makebox(0,0)[l]{$S\langle C,A\rangle$}}
\put(40,40){\makebox(0,0)[tr]{\scriptsize$T\langle \mj_{A},f\rangle$}}
\put(40,70){\makebox(0,0)[br]{\scriptsize$T\langle f,\mj_{C}\rangle$}}
\put(160,40){\makebox(0,0)[tl]{\scriptsize$S\langle f,\mj_{A}\rangle$}}
\put(160,70){\makebox(0,0)[bl]{\scriptsize$S\langle \mj_{C},f\rangle$}}
\put(100,18){\makebox(0,0)[t]{\scriptsize$\alpha(A)$}}
\put(100,94){\makebox(0,0)[b]{\scriptsize$\alpha(C)$}}
\end{picture}
\end{center}
where $\langle u,v \rangle$ denotes the tuple of arguments
with $u$ in positive and $v$ in negative argument places.

Let now $\Gamma$ be a graph with $k$ ($k>1$) components, and
let
\[\alpha=\{\alpha(A_1,\ldots ,A_k):T(A_{\pi(x_1)},\ldots ,A_{\pi(x_m)})
\str S(A_{\pi(y_1)},\ldots ,A_{\pi(y_n)})\; | \; A_1,\ldots A_k \in
Ob({\cal A})\}\]
be a family of morphisms from $\cal A$ indexed by the set of $k$-tuples
of objects from $\cal A$. Then we say that $\alpha$ is a g-dinatural
transformation from $T$ to $S$ with the graph $\Gamma$, if for
every $k-1$-tuple $(A_1,\ldots ,A_{i-1},A_{i+1},\ldots ,A_k)$ of
objects from $\cal A$, the subset
\[\alpha_{A_1,\ldots ,A_{i-1},A_{i+1},\ldots ,A_k}=
\{\alpha(A_1,\ldots ,A_{i-1},A,A_{i+1},\ldots ,A_k)\; | \; A\in Ob({\cal A})\}
\]
of $\alpha$ is g-dinatural with the graph $\Gamma_i$. (All the argument
places that are not in $\Gamma_i$ are parametrized in this case.)
This means that a transformation is g-dinatural iff it is g-dinatural
in each of its components, or roughly speaking, g-dinaturality is
defined componentwise.

\exmp Let $\cal C$ be a cartesian closed category, and let
$T:{\cal C}\times{\cal C}^{op}\times{\cal C}\str{\cal C}$ and
$R:{\cal C}\times{\cal C}^{op}\times{\cal C}\times{\cal C}\str{\cal C}$ be two
functors defined on objects and morphisms of $\cal C$ by the formulae
\[T(x_1,x_2,x_3)=x_1\times(x_2\str x_3)\quad{\mbox{\rm and}}\quad
R(z_1,z_2,z_3,z_4)=(z_1\times(z_2\str z_3))\times z_4.\]
Let $\Gamma$ be the graph
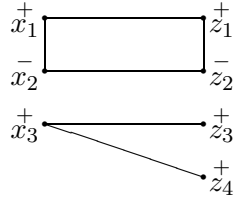
\begin{center}
\begin{picture}(100,100)
\put(20,40){\circle*{2}}
\put(18,40){\makebox(0,0)[r]{$\stackrel{+}{x_{3}}$}}
\put(20,60){\circle*{2}}
\put(18,60){\makebox(0,0)[r]{$\stackrel{-}{x_{2}}$}}
\put(20,80){\circle*{2}}
\put(18,80){\makebox(0,0)[r]{$\stackrel{+}{x_{1}}$}}
\put(80,20){\circle*{2}}
\put(82,20){\makebox(0,0)[l]{$\stackrel{+}{z_{4}}$}}
\put(80,40){\circle*{2}}
\put(82,40){\makebox(0,0)[l]{$\stackrel{+}{z_{3}}$}}
\put(80,60){\circle*{2}}
\put(82,60){\makebox(0,0)[l]{$\stackrel{-}{z_{2}}$}}
\put(80,80){\circle*{2}}
\put(82,80){\makebox(0,0)[l]{$\stackrel{+}{z_{1}}$}}
\put(20,40){\line(1,0){60}}
\put(20,40){\line(3,-1){60}}
\put(20,60){\line(1,0){60}}
\put(20,60){\line(0,1){20}}
\put(20,80){\line(1,0){60}}
\put(80,60){\line(0,1){20}}
\end{picture}
\end{center}
and let $\alpha$ be the following transformation
\[\{\alpha(A,B)=(\mj_{A\times(A\str B)}\times\mep_{A,B})\mw_{A\times(A\str B)}:
A\times(A\str B) \vdash (A\times(A\str B))\times B \; | \; A,B\in Ob({\cal C})
\}.\]
Then $\alpha:T\gdna R$, because for all $A$, $B$, $C$ and $f:A\vdash C$
from $\cal C$
the following two diagrams commute:

\begin{center}
\begin{picture}(200,110)
\put(80,20){\line(1,0){40}}
\put(80,90){\line(1,0){40}}
\put(30,50){\line(3,-2){30}}
\put(140,80){\line(3,-2){30}}
\put(30,60){\line(3,2){30}}
\put(140,30){\line(3,2){30}}
\put(80,17){\line(0,1){6}}
\put(80,87){\line(0,1){6}}
\put(30,47){\line(0,1){6}}
\put(30,57){\line(0,1){6}}
\put(140,27){\line(0,1){6}}
\put(140,77){\line(0,1){6}}
\put(78,20){\makebox(0,0)[r]{$A\times(A{\str}B)$}}
\put(78,90){\makebox(0,0)[r]{$C\times(C{\str}B)$}}
\put(122,20){\makebox(0,0)[l]{$(A\times(A{\str}B)){\times}B$}}
\put(122,90){\makebox(0,0)[l]{$(C\times(C{\str}B)){\times}B$}}
\put(25,55){\makebox(0,0)[r]{$A\times(C{\str}B)$}}
\put(175,55){\makebox(0,0)[l]{$(C\str(A{\str}B)){\times}B$}}
\put(40,40){\makebox(0,0)[tr]{\scriptsize$\mj_A\times(f{\str}\mj_{B})$}}
\put(40,70){\makebox(0,0)[br]{\scriptsize$f\times(\mj_{C}{\str}\mj_{B})$}}
\put(160,40){\makebox(0,0)[tl]{\scriptsize$(f\times(\mj_{A}{\str}\mj_{B}))
{\times}\mj_{B}$}}
\put(160,70){\makebox(0,0)[bl]{\scriptsize$(\mj_C\times(f{\str}\mj_{B}))
{\times}\mj_{B}$}}
\put(100,18){\makebox(0,0)[t]{\scriptsize$\alpha(A,B)$}}
\put(100,93){\makebox(0,0)[b]{\scriptsize$\alpha(C,B)$}}
\end{picture}
\end{center}
\begin{center}
\begin{picture}(200,110)
\put(80,20){\line(1,0){40}}
\put(80,90){\line(1,0){40}}
\put(30,50){\line(3,-2){30}}
\put(140,80){\line(3,-2){30}}
\put(30,60){\line(3,2){30}}
\put(140,30){\line(3,2){30}}
\put(80,17){\line(0,1){6}}
\put(80,87){\line(0,1){6}}
\put(30,47){\line(0,1){6}}
\put(30,57){\line(0,1){6}}
\put(140,27){\line(0,1){6}}
\put(140,77){\line(0,1){6}}
\put(78,20){\makebox(0,0)[r]{$B\times(B{\str}A)$}}
\put(78,90){\makebox(0,0)[r]{$B\times(B{\str}C)$}}
\put(122,20){\makebox(0,0)[l]{$(B\times(B{\str}A)){\times}A$}}
\put(122,90){\makebox(0,0)[l]{$(B\times(B{\str}C)){\times}C$}}
\put(25,55){\makebox(0,0)[r]{$B\times(B{\str}A)$}}
\put(175,55){\makebox(0,0)[l]{$(B\times(B{\str}C)){\times}C$}}
\put(40,40){\makebox(0,0)[tr]{\scriptsize$\mj_B\times(\mj_{B}{\str}\mj_{A})$}}
\put(40,70){\makebox(0,0)[br]{\scriptsize$\mj_B\times(\mj_{B}{\str}f)$}}
\put(160,40){\makebox(0,0)[tl]{\scriptsize$(\mj_B\times(\mj_{B}{\str}f))
{\times}f$}}
\put(160,70){\makebox(0,0)[bl]{\scriptsize$(\mj_B\times(\mj_{B}{\str}\mj_{C}))
{\times}\mj_{C}$}}
\put(100,18){\makebox(0,0)[t]{\scriptsize$\alpha(B,A)$}}
\put(100,93){\makebox(0,0)[b]{\scriptsize$\alpha(B,C)$}}
\end{picture}
\end{center}

\exmp Let $\beta$ be the following transformation
\[
\{\beta(A,B)=((\mk_A\times \mj_{A\str B})\times\mj_B)(\mj_{A\times(A\str B)}
\times\mep_{A,B})\mw_{A\times(A\str B)}(\mj_A\times(\mk_A\str\mj_B))\; | \;
A,B\in{\cal C}\}\]
between the functors $T:{\cal C}\times{\cal C}\str\cal C$ and
$S:{\cal C}^{op}\times{\cal C}\times{\cal C}\str\cal C$ that are defined
by the terms $x_1\times(\ri\str x_2)$ and
$(\ri\times(y_1\str y_2))\times y_3$ respectively,
for some cartesian closed category $\cal C$.
Then we can show that $\beta$ is g-dinatural with the graph:
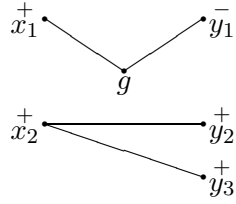
\begin{center}
\begin{picture}(100,100)
\put(20,40){\circle*{2}}
\put(18,40){\makebox(0,0)[r]{$\stackrel{+}{x_{2}}$}}
\put(20,80){\circle*{2}}
\put(18,80){\makebox(0,0)[r]{$\stackrel{+}{x_{1}}$}}
\put(80,20){\circle*{2}}
\put(82,20){\makebox(0,0)[l]{$\stackrel{+}{y_{3}}$}}
\put(80,40){\circle*{2}}
\put(82,40){\makebox(0,0)[l]{$\stackrel{+}{y_{2}}$}}
\put(50,60){\circle*{2}}
\put(50,58){\makebox(0,0)[t]{$g$}}
\put(80,80){\circle*{2}}
\put(82,80){\makebox(0,0)[l]{$\stackrel{-}{y_{1}}$}}
\put(20,40){\line(1,0){60}}
\put(20,40){\line(3,-1){60}}
\put(20,80){\line(3,-2){30}}
\put(80,80){\line(-3,-2){30}}
\end{picture}
\end{center}

\vspace{.2cm}
It is obvious how the notion of g-dinaturality extends the traditional notion
of dinaturality given in \cite{DS70}.
All that one has to do in order to show
that a g-dinatural transformation is already dinatural is to collapse all the
argument places of the same sign from a component into
one argument place.
The main purpose of this extension is to
give an answer to the question: ``When is the composition of two dinatural
transformations dinatural?'' The rest of this section
is devoted to this problem.

Let $\Phi=(V_\Phi,l_\Phi,E_\Phi)$ where $V_\Phi=\{x_1,\ldots,x_m,y_1,\ldots,
y_n\}\cup G_\Phi$ and
$\Psi=(V_\Psi,l_\Psi,E_\Psi)$ where $V_\Psi=\{y_1,\ldots,y_n,z_1,\ldots,
z_p\}\cup G_\Psi$ be two graphs with $k_\Phi$ and $k_\Psi$ components,
respectively, such that $l_\Phi$ and $l_\Psi$ coincide on
$\{y_1,\ldots,y_n\}$ and that $G_\Phi\cap G_\Psi = \emptyset$.
Let
\[\begin{array}{l}
T:{\cal A}^{l_\Phi(x_1)}\times {\cal A}^{l_\Phi(x_2)}\times\ldots\times
{\cal A}^{l_\Phi(x_m)}\str{\cal A}\\[.2cm]
S:{\cal A}^{l_\Phi(y_1)}\times {\cal A}^{l_\Phi(y_2)}\times\ldots\times
{\cal A}^{l_\Psi(y_n)}\str{\cal A}\\[.2cm]
R:{\cal A}^{l_\Psi(z_1)}\times {\cal A}^{l_\Psi(z_2)}\times\ldots\times
{\cal A}^{l_\Psi(z_p)}\str{\cal A}
\end{array}\]
be three functors, and let $\alpha$ and $\beta$ be two g-dinatural
transformations
\[\begin{array}{l}
\alpha=\{\alpha(A_1,\ldots,A_{k_\Phi})\; | \; A_1,\ldots,A_{k_\Phi}\in
Ob({\cal A})\}:T\fdna S\\[.2cm]
\beta=\{\beta(B_1,\ldots,B_{k_\Psi})\; | \; B_1,\ldots,B_{k_\Psi}\in
Ob({\cal A})\}:S\pdna R
\end{array}\]
By the {\em amalgamation} of $\Phi$ and $\Psi$ we mean the couple
$(V_\Phi \cup V_\Psi\: ,\: E_\Phi + E_\Psi)$ denoted by $\Phi+\Psi$.
(Note that $\Phi+\Psi$ is not a graph in the sense of the definition above,
but we may define its components analogously.)

Let the amalgamation $\Phi+\Psi$ have one component; then we define the
graph $\Psi\Phi=(V,l,E)$, i.e. the {\em composition of the graphs}
$\Phi$ and $\Psi$,  in the following manner:
\\[.1cm]
--  if all $x$'s are of the same sign in $\Phi$, which is opposite to
the sign of all $z$'s in $\Psi$
(this includes the cases when $m=0$ or $p=0$), then $V=\{x_1,\ldots,x_m\}
\cup\{g\}$ and $E=\{\{x_i,g\}\; | \; 1\leq i \leq m\} \cup
\{\{z_j,g\}\; | \; 1\leq j \leq p\}$,
\\[.1cm]
-- otherwise, $G$ is empty, $V=\{x_1,\ldots,x_m,z_1,\ldots,z_p\}$ and
$E=\{\{x_i,x_j\}\; | \; l_\Phi(x_i)=-l_\Phi(x_j)\}\cup
\{\{z_i,z_j\}\; | \; l_\Psi(z_i)=-l_\Psi(z_j)\}\cup
\{\{x_i,z_j\}\; | \; l_\Phi(x_i)=l_\Psi(z_j)\}$.
\\[.2cm]
In both cases, the function $l$ is defined so that its restrictions to
$\{x_1,\ldots,x_m\}$ and $\{z_1,\ldots,z_p\}$ are $l_\Phi$ and $l_\Psi$
respectively.

In the case of more than one component in $\Phi+\Psi$, we
proceed analogously for each of them to construct a component
of the graph $\Psi\Phi$. Since the notion of g-dinaturality is defined
componentwise, from now on we consider just the case when $\Phi+\Psi$,
and therefore $\Psi\Phi$, has a single component.

Now we define the composition $\beta\alpha$ to be the transformation
\[\{\beta(A^{k_\Psi})\alpha(A^{k_\Phi})\; | \; A\in Ob({\cal A})\}.\]
Our question is: ``Is it a g-dinatural transformation with the graph
$\Psi\Phi$?''

\exmp Let $\cal C$, $T$, $R$ be as in Example \arabic{section}.2 and let
$S:{\cal A}\times{\cal A}^{op}\times{\cal A}\times{\cal A}\times
{\cal A}^{op}\times{\cal A}\str{\cal A}$ be defined on objects and morphisms
by the formula
\[S(y_1,y_2,y_3,y_4,y_5,y_6)=(y_1\times(y_2\str y_3))\times(y_4\times
(y_5\str y_6)).\]
Let $\Phi$ be the thin graph and $\Psi$ the thick graph whose amalgamation
$\Phi+\Psi$ is given by the diagram
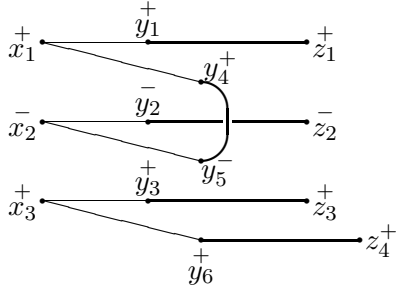
\begin{center}
\begin{picture}(160,120)
\put(20,40){\circle*{2}}
\put(18,40){\makebox(0,0)[r]{$\stackrel{+}{x_{3}}$}}
\put(20,70){\circle*{2}}
\put(18,70){\makebox(0,0)[r]{$\stackrel{-}{x_{2}}$}}
\put(20,100){\circle*{2}}
\put(18,100){\makebox(0,0)[r]{$\stackrel{+}{x_{1}}$}}
\put(140,25){\circle*{2}}
\put(142,25){\makebox(0,0)[l]{$z^+_4$}}
\put(120,40){\circle*{2}}
\put(122,40){\makebox(0,0)[l]{$\stackrel{+}{z_{3}}$}}
\put(120,70){\circle*{2}}
\put(122,70){\makebox(0,0)[l]{$\stackrel{-}{z_{2}}$}}
\put(120,100){\circle*{2}}
\put(122,100){\makebox(0,0)[l]{$\stackrel{+}{z_{1}}$}}
\put(60,40){\circle*{2}}
\put(60,44){\makebox(0,0)[b]{$\stackrel{+}{y_{3}}$}}
\put(60,70){\circle*{2}}
\put(60,74){\makebox(0,0)[b]{$\stackrel{-}{y_{2}}$}}
\put(60,100){\circle*{2}}
\put(60,104){\makebox(0,0)[b]{$\stackrel{+}{y_{1}}$}}
\put(80,25){\circle*{2}}
\put(80,23){\makebox(0,0)[t]{$\stackrel{+}{y_{6}}$}}
\put(80,55){\circle*{2}}
\put(80,58){\makebox(0,0)[tl]{$y^-_5$}}
\put(80,85){\circle*{2}}
\put(81,87){\makebox(0,0)[bl]{$y^+_4$}}
{\linethickness{.01pt}\put(20,40){\line(1,0){40}}
\put(20,40){\line(4,-1){60}}
\put(20,70){\line(1,0){40}}
\put(20,70){\line(4,-1){60}}
\put(20,100){\line(1,0){40}}
\put(20,100){\line(4,-1){60}}}
{\thicklines
\put(80,25){\line(1,0){60}}
\put(60,40){\line(1,0){60}}
\put(60,70){\line(1,0){28}}
\put(92,70){\line(1,0){28}}
\put(60,100){\line(1,0){60}}
\put(80,70){\oval(20,30)[r]}
}
\end{picture}
\end{center}
Let $\beta$ and $\gamma$ be the transformations
\[\begin{array}{l}
\{\beta(A,B,C)=\mw_{A\times(B\str C)}\; | \; A,B,C\in Ob({\cal C})\}\\[.2cm]
\{\gamma(A,B,C,D,Z)=\mj_{A\times(B\str C)}\times \mep_{D,Z} \; | \;
A,B,C,D,Z\in Ob({\cal C})\}
\end{array}\]
Then it is easy to check that $\beta:T\fdna S$, $\gamma:S\pdna R$ and that
$\Psi\Phi=\Gamma$ and $\gamma\beta=\alpha$ for $\Gamma$ and $\alpha$ from
Example \arabic{section}.2.

\vspace{.2cm}
One may be tempted by these examples to conclude that
the composition of g-dinatural transformations is
always g-dinatural, as it is the case with natural transformations. This will
be proven wrong. However, the category in question may have strong influence
on g-dinaturality of the composition of g-dinatural transformations, but
we will neglect this possible influence and rely only on the geometry
of the underlying graphs. An approach that treats properties intrinsic to
a category that are sufficient for dinaturality of a composition of
transformations is given in \cite{BFSS90}.

The next example, although tedious, may serve as a good introduction to
what follows.

\exmp Let
$T:{\cal A}\times{\cal A}^{op}\times{\cal A}^2\str{\cal A}$, 
$S:{\cal A}^2\times{({\cal A}^{op})}^2\times{\cal A}\times
{({\cal A}^{op})}^2\times{\cal A}\str{\cal A}$ and $R:{\cal A}\str{\cal A}$
be three functors and $\alpha:T\fdna S$ and $\beta:S\pdna R$ two
g-dinatural transformations such that the amalgamation $\Phi+\Psi$ ($\Psi$ is
bold) is given by the following diagram
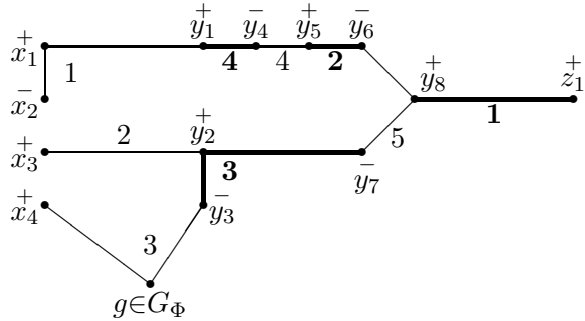
\begin{center}
\begin{picture}(240,130)
\put(20,50){\circle*{3}}
\put(18,50){\makebox(0,0)[r]{$\stackrel{+}{x_{4}}$}}
\put(20,70){\circle*{3}}
\put(18,70){\makebox(0,0)[r]{$\stackrel{+}{x_{3}}$}}
\put(20,90){\circle*{3}}
\put(18,90){\makebox(0,0)[r]{$\stackrel{-}{x_{2}}$}}
\put(20,110){\circle*{3}}
\put(18,110){\makebox(0,0)[r]{$\stackrel{+}{x_{1}}$}}
\put(80,70){\circle*{3}}
\put(80,75){\makebox(0,0)[b]{$\stackrel{+}{y_{2}}$}}
\put(80,110){\circle*{3}}
\put(80,115){\makebox(0,0)[b]{$\stackrel{+}{y_{1}}$}}
\put(80,50){\circle*{3}}
\put(82,50){\makebox(0,0)[l]{$\stackrel{-}{y_{3}}$}}
\put(100,110){\circle*{3}}
\put(100,115){\makebox(0,0)[b]{$\stackrel{-}{y_{4}}$}}
\put(120,110){\circle*{3}}
\put(120,115){\makebox(0,0)[b]{$\stackrel{+}{y_{5}}$}}
\put(140,70){\circle*{3}}
\put(142,69){\makebox(0,0)[t]{$\stackrel{-}{y_{7}}$}}
\put(140,110){\circle*{3}}
\put(140,115){\makebox(0,0)[b]{$\stackrel{-}{y_{6}}$}}
\put(160,90){\circle*{3}}
\put(162,95){\makebox(0,0)[bl]{$\stackrel{+}{y_{8}}$}}
\put(220,90){\circle*{3}}
\put(220,95){\makebox(0,0)[b]{$\stackrel{+}{z_{1}}$}}
\put(60,20){\circle*{3}}
\put(60,17){\makebox(0,0)[t]{$g{\in}G_{\Phi}$}}
\put(30,100){\makebox(0,0){$1$}}
\put(50,73){\makebox(0,0)[b]{$2$}}
\put(60,35){\makebox(0,0){$3$}}
\put(110,108){\makebox(0,0)[t]{$4$}}
\put(151,79){\makebox(0,0)[tl]{$5$}}
\put(190,88){\makebox(0,0)[t]{\bf 1}}
\put(130,108){\makebox(0,0)[t]{\bf 2}}
\put(90,63){\makebox(0,0){\bf 3}}
\put(90,108){\makebox(0,0)[t]{\bf 4}}
\put(20,50){\line(4,-3){40}}
\put(20,70){\line(1,0){60}}
\put(20,90){\line(0,1){20}}
\put(20,110){\line(1,0){60}}
\put(60,20){\line(2,3){20}}
\put(100,110){\line(1,0){20}}
\put(140,110){\line(1,-1){20}}
\put(140,70){\line(1,1){20}}
{\linethickness{1.5pt}
\put(80,70){\line(0,-1){20}}
\put(80,70){\line(1,0){60}}
\put(80,110){\line(1,0){20}}
\put(120,110){\line(1,0){20}}
\put(160,90){\line(1,0){60}}
}
\end{picture}
\end{center}
where the components of $\Phi$ and $\Psi$ are enumerated by suitable numerals.
The composition of $\Phi$ and $\Psi$ is given by the diagram

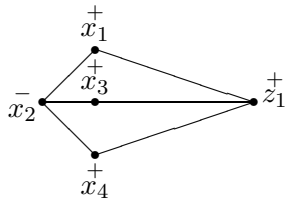
\begin{center}
\begin{picture}(120,80)
\put(20,40){\circle*{3}}
\put(18,40){\makebox(0,0)[r]{$\stackrel{-}{x_{2}}$}}
\put(40,20){\circle*{3}}
\put(40,18){\makebox(0,0)[t]{$\stackrel{+}{x_{4}}$}}
\put(40,40){\circle*{3}}
\put(40,45){\makebox(0,0)[b]{$\stackrel{+}{x_{3}}$}}
\put(40,60){\circle*{3}}
\put(40,65){\makebox(0,0)[b]{$\stackrel{+}{x_{1}}$}}
\put(100,40){\circle*{3}}
\put(103,40){\makebox(0,0)[bl]{$\stackrel{+}{z_{1}}$}}
\put(20,40){\line(1,-1){20}}
\put(20,40){\line(1,0){80}}
\put(20,40){\line(1,1){20}}
\put(40,20){\line(3,1){60}}
\put(40,60){\line(3,-1){60}}
\end{picture}
\end{center}
and $\beta\alpha$ is g-dinatural with this graph if the following
equation
\[R(\mj_C)\beta\alpha(C)T(f,\mj_C,f^2)=R(f)\beta\alpha(A)T(\mj_A,f,{\mj_A}^2)
\]
holds in $\cal A$ for every $A$, $C$ and $f:A \str C$ from this category.
We prove this by ``travelling'' along the amalgamation $\Phi+\Psi$,
relying on the definition of $\beta\alpha$, on the functoriality of $T$, $S$
and
$R$ and on the g-dinaturality of $\alpha$ and $\beta$. We hope 
the reader won't be scared with the following a rather long proof
in which $(\beta\alpha)$ means reference to the definition of $\beta\alpha$,
$(T)$ means reference to functoriality of $T$, $(\alpha 3)$ means
reference to g-dinaturality of
$\alpha$ in the third component of $\Phi$, etc.
\[\begin{array}{l}
R(\mj_C)\beta\alpha(C)T(f,\mj_C,f^2)\\[0.15cm]
=R(\mj_{C})\beta(C^{4})\alpha(C^{5})T(f,\mj_{C},f^{2})\quad (\beta\alpha)\\[0.15cm]
=R(\mj_{C})\beta(C^{4})\alpha(C^{5})T(f,\mj_{C}^{3})
T(\mj_{A},\mj_{C},f^{2})\quad (T)\\[0.15cm]
=R(\mj_{C})\beta(C^{4})S(f,\mj_{C}^{7})\alpha(A,C^{4})
T(\mj_{A},f,\mj_{C}^{2})T(\mj_{A},\mj_{C},f^{2})\quad (\alpha1)\\[0.15cm]
=R(\mj_{C})\beta(C^{3},A)S(\mj_{A},\mj_{C}^{2},f,\mj_{C}^{4})
\alpha(A,C^{4})T(\mj_{A},f^{3})\quad (\beta4),(T)\\[0.15cm]
=R(\mj_{C})\beta(C^{3},A)S(\mj_{A},\mj_{C}^{2},\mj_{A},f,\mj_{C}^{3})
\alpha(A,C^{2},A,C)T(\mj_{A},f^{3})\quad (\alpha4)\\[0.15cm]
=R(\mj_{C})\beta(C,A,C,A)S(\mj_{A},\mj_{C}^{2},\mj_{A}^{2},f,\mj_{C}^{2})
\alpha(A,C^{2},A,C)T(\mj_{A},f^{3})\quad (\beta2)\\[0.15cm]
=R(\mj_{C})\beta(C,A,C,A)S(\mj_{A},\mj_{C}^{2},\mj_{A}^{2},f,\mj_{C}^{2})
S(\mj_A,\mj_C^2,\mj_A^2,\mj_C^3)\alpha(A,C^{2},A,C)T(\mj_{A}^{2},f,\mj_{C})\\
{\mbox{\hspace{7em}}}T(\mj_{A},f,\mj_{A},f)\quad(T),(S)\\[0.15cm]
=R(\mj_{C})\beta(C,A,C,A)S(\mj_{A},\mj_{C}^{2},\mj_{A}^{2},f,\mj_{C}^{2})
S(\mj_{A},f,\mj_{C},\mj_{A}^{2},\mj_{C}^{3})\alpha(A^{2},C,A,C)\\
{\mbox{\hspace{7em}}}
T(\mj_{A}^{3},\mj_{C})T(\mj_{A},f,\mj_{A},f)\quad (\alpha2)\\[0.15cm]
=R(\mj_{C})\beta(C,A,C,A)S(\mj_{A},f,\mj_{C},\mj_{A}^{3},\mj_{C}^{2})
S(\mj_{A}^{2},\mj_{C},\mj_{A}^{2},f,\mj_{C}^{2})\alpha(A^{2},C,A,C)\\
{\mbox{\hspace{7em}}}
T(\mj_{A},f,\mj_{A},f)\quad (S),(T)\\[0.15cm]
=R(\mj_{C})\beta(C,A^{3})S(\mj_{A}^{2},f,\mj_{A}^{3},f,\mj_{C})
S(\mj_{A}^{2},\mj_{C},\mj_{A}^{2},f,\mj_{C}^{2})
\alpha(A^{2},C,A,C)T(\mj_{A},f,\mj_{A},f)\quad (\beta3)\\[0.15cm]
=R(\mj_{C})\beta(C,A^{3})S(\mj_{A}^{5},f^{2},\mj_{C})
S(\mj_{A}^{2},f,\mj_{A}^{2},\mj_{C}^{3})\alpha(A^{2},C,A,C)
T(\mj_{A}^{3},f)T(\mj_{A},f,\mj_{A}^{2})\quad (T),(S)\\[0.15cm]
=R(\mj_{C})\beta(C,A^{3})S(\mj_{A}^{5},f^{2},\mj_{C})
S(\mj_{A}^{5},\mj_{C}^{3})\alpha(A^{4},C)T(\mj_{A}^{4})
T(\mj_{A},f,\mj_{A}^{2})\quad (\alpha3)\\[0.15cm]
=R(\mj_{C})\beta(C,A^{3})S(\mj_{A}^{5},f^{2},\mj_{C})
\alpha(A^{4},C)T(\mj_{A}^{4})T(\mj_{A},f,\mj_{A}^{2})\quad (S)\\[0.15cm]
=R(\mj_{C})\beta(C,A^{3})S(\mj_{A}^{7},f)
\alpha(A^{5})T(\mj_{A},f,\mj_{A}^{2})\quad (\alpha5),(T)\\[0.15cm]
=R(f)\beta(A^{4})\alpha(A^{5})T(\mj_{A},f,\mj_{A}^{2})\quad (\beta1),(S)
\\[0.15cm]
=R(f)\beta\alpha(A)T(\mj_A,f,{\mj_A}^2)\quad (\beta\alpha)
\end{array}\]
This example strengthens the impression that g-dinatural transformations
give a g-dinatural transformation in the composition, but could we repeat
the above procedure with transformations whose amalgamation of graphs is
given below?
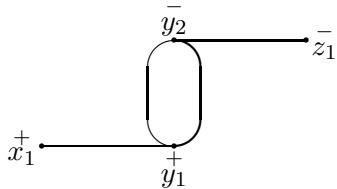
\begin{center}
\begin{picture}(140,80)
\put(20,20){\circle*{2}}
\put(18,20){\makebox(0,0)[r]{$\stackrel{+}{x_{1}}$}}
\put(70,20){\circle*{2}}
\put(70,19){\makebox(0,0)[t]{$\stackrel{+}{y_{1}}$}}
\put(70,60){\circle*{2}}
\put(70,64){\makebox(0,0)[b]{$\stackrel{-}{y_{2}}$}}
\put(120,60){\circle*{2}}
\put(122,60){\makebox(0,0)[l]{$\stackrel{-}{z_{1}}$}}
\put(20,20){\line(1,0){50}}
\put(70,40){\oval(20,40)[l]}
{\thicklines
\put(70,60){\line(1,0){50}}
\put(70,40){\oval(20,40)[r]}
}
\end{picture}
\end{center}
Simply, without any further assumptions on the category in question, we
can't move along this amalgamation at all.

We shall now examine properties of an amalgamation $\Phi+\Psi$ which
guarantee that the composition of $\alpha:T\fdna S$ and $\beta:S\pdna R$
is g-dinatural. For these purposes let $\Phi$ and $\Psi$ be as in the
definition of amalgamation, and let $\Phi+\Psi$ have one component.
We say that $\Phi+\Psi$ {\em provides g-dinaturality} if for every
category $\cal A$, for every triple of functors 
$F:{\cal A}^{l_\Phi(x_1)}\times\ldots\times{\cal A}^{l_\Phi(x_m)}
\str{\cal A}$,
$G:{\cal A}^{l_\Phi(y_1)}\times\ldots\times{\cal A}^{l_\Psi(y_n)}
\str{\cal A}$ and
$H:{\cal A}^{l_\Psi(z_1)}\times\ldots\times{\cal A}^{l_\Psi(z_p)}
\str{\cal A}$,
and for every pair $\gamma:F\fdna G$ and $\delta:G\pdna H$ of g-dinatural
transformations, the composition $\delta\gamma$ is g-dinatural from
$F$ to $H$ with the graph $\Psi\Phi$. Let $P(\Phi,\Psi)$
denote the property that $\Phi+\Psi$ provides g-dinaturality.
To make easier the proof of the main result of this section, we introduce
an alternative characterization of $P(\Phi,\Psi)$. In the style of
\cite{D99} we introduce a free categorial object that will serve as a
template for g-dinaturality.

Let ${\cal K}_{\Phi,\Psi}$ be the category of structured categories
$({\cal A},F,G,H,\gamma,\delta)$ for $\cal A$, $F$, $G$, $H$, $\gamma$,
$\delta$ as above. The morphisms of ${\cal K}_{\Phi,\Psi}$ are
structure-preserving functors between these categories. The category
${\cal K}_{\Phi,\Psi}$ has an equational presentation, as we shall
see later; hence, there exists a free object of this category generated
by the arrow
\[ A \stackrel{\mf}{\longrightarrow} C.\]
Denote this object by $({\cal D},T,S,R,\alpha,\beta)$. Its explicit
construction will be given soon. The following lemma gives an
alternative definition of $P(\Phi,\Psi)$.

\lema {\em The amalgamation} $\Phi+\Psi$ {\em provides g-dinaturality
iff the following diagram}
\begin{center}
\begin{picture}(220,110)
\put(80,20){\vector(1,0){60}}
\put(80,90){\vector(1,0){60}}
\put(30,50){\vector(3,-2){30}}
\put(160,80){\vector(3,-2){30}}
\put(30,60){\vector(3,2){30}}
\put(160,30){\vector(3,2){30}}
\put(78,20){\makebox(0,0)[r]{$T\langle A,A\rangle$}}
\put(78,90){\makebox(0,0)[r]{$T\langle C,C\rangle$}}
\put(142,20){\makebox(0,0)[l]{$R\langle A,A\rangle$}}
\put(142,90){\makebox(0,0)[l]{$R\langle C,C\rangle$}}
\put(25,55){\makebox(0,0)[r]{$T\langle A,C\rangle$}}
\put(195,55){\makebox(0,0)[l]{$R\langle C,A\rangle$}}
\put(40,40){\makebox(0,0)[tr]{\scriptsize$T\langle \mj_{A},\mf\rangle$}}
\put(40,70){\makebox(0,0)[br]{\scriptsize$T\langle \mf,\mj_{C}\rangle$}}
\put(180,40){\makebox(0,0)[tl]{\scriptsize$R\langle \mf,\mj_{A}\rangle$}}
\put(180,70){\makebox(0,0)[bl]{\scriptsize$R\langle \mj_{C},\mf\rangle$}}
\put(110,18){\makebox(0,0)[t]{\scriptsize$\beta(A^{k_{\beta}})
\alpha(A^{k_{\alpha}})$}}
\put(113,96){\makebox(0,0)[b]{\scriptsize$\beta(C^{k_{\beta}})
\alpha(C^{k_{\alpha}})$}}
\end{picture}
\end{center}
{\em commutes in} $\cal D$, {\em where} $\mf:A\str C$ {\em is the generator of}
$\cal D$.

\dkz The ``only if'' part of the lemma follows from the definitions of
g-dinaturality and of $P(\Phi,\Psi)$. For the ``if'' part we rely on the
universal property of the category $\cal D$.
\qed

The category $\cal D$ can be built up from syntactical material
in the following manner.
The {\em objects} of $\cal D$ are freely generated over the set $\{A,C\}$
by the $m$-ary operation $T$, the $n$-ary operation $S$ and the $p$-ary operation
$R$. We use the schematic letters $X$, $Y$ and $Z$, possibly with indices,
for elements of $Ob({\cal D})$.
The {\em primitive morphism terms} of $\cal D$ are
\[\begin{array}{l}
\mf:A\str C,\quad\quad\quad \mj_X:X\str X,\\[.2cm]
\alpha(Y_1,\ldots,Y_{k_\Phi}):T(Y_{\pi(x_1)},\ldots,Y_{\pi(x_m)})\str
S(Y_{\pi(y_1)},\ldots,Y_{\pi(y_n)}),
\\[.2cm]
\beta(Z_1,\ldots,Z_{k_\Psi}):S(Z_{\pi'(y_1)},\ldots,Z_{\pi'(y_m)})\str
R(Z_{\pi'(z_1)},\ldots,Z_{\pi'(z_p)}),
\end{array}\]
for all objects $X,Y_1,\ldots,Y_{k_\Phi},Z_1,\ldots,Z_{k_\Psi}$, where $\pi$
and $\pi'$ are component classifiers for $\Phi$ and $\Psi$, respectively.

In the following definitions and equations let $F$ range over the set
$\{T, S, R\}$, and let $k$, depending on $F$, be the
variable for $m$, $n$ or $p$
respectively.

{\em Morphism terms} of $\cal D$ are defined inductively as follows:
\begin{enumerate}
\item primitive morphism terms are morphism terms,
\item if $g:X\str Y$ and $h:Y\str Z$ are morphism terms, then
$hg:X\str Z$ is a morphism term,
\item if $\{t_i:X_i\str Y_i \; | \; 1\leq i\leq k\; {\mbox{\rm and the $i$-th
argument place of $F$ is positive}}\}$ and
$\{t_j:Y_j\str X_j \; | \; 1\leq j\leq k,\; {\mbox{\rm and the $j$-th
argument place of $F$ is negative}}\}$
are two sets of morphism terms, then
$F(t_1,\ldots,t_m):F(X_1,\ldots,X_m) \str F(Y_1,\ldots,Y_m)$ is a morphism
term.
\end{enumerate}
For morphism terms we use the schematic letters $g$, $h$, $t$, possibly
primed and with indices,
and $\equiv$ is used for identity of terms.
{\em Morphisms} of $\cal D$ are the equivalence classes of morphism terms
modulo congruence generated by the following schematic equations.
\\[.2cm]
{\em Categorial equations}
\\[.1cm]
{\makebox(50,0)[bl]{$(cat1)$}}$g\mj_{X}=g=\mj_Y\: g$.
\\[.1cm]
{\makebox(50,0)[bl]{$(cat2)$}}$t(hg)=(th)g$.
\\[.2cm]
{\em Functorial equations}
\\[.1cm]
{\makebox(50,0)[bl]{$(F)$}}For morphism terms $g_1,\ldots,g_k,h_1,\ldots,
h_k,t_1,\ldots,t_k$ such that for every $1\leq i\leq k$
\[t_i\equiv\left\{
\begin{array}{ll}
h_ig_i & {\mbox{\rm ; if the $i$-th argument place of $F$ is positive}}\\
g_ih_i & {\mbox{\rm ; if the $i$-th argument place of $F$ is negative}}
\end{array},\right.\]
{\makebox(50,0)[bl]{}}$F(h_1,\ldots,h_k)F(g_1,\ldots,g_k)=F(t_1,\ldots,t_k)$.
\\[.1cm]
{\makebox(50,0)[bl]{$(F1)$}}$F(\mj_{X_1},\ldots,\mj_{X_k})=
\mj_{F(X_1,\ldots,X_k)}$.
\\[.2cm]
{\em G-dinatural equations}
\\[.1cm]
{\makebox(50,0)[bl]{$(\alpha)$}}For $1\leq i\leq k_\Phi$ and
morphism terms $t:X\str Y$, $g_1,\ldots,g_m$,
$g_1',\ldots,g_m'$, $h_1,\ldots,h_n$, $h_1',\ldots,h_n'$ such that for
$1\leq j\leq m$
\[g_j\equiv\left\{
\begin{array}{ll}
\mj_{Z_q} & ;\; x_j\in\Phi_q\neq\Phi_i \\
t         & ;\; x_j\in\Phi^+_i \\
\mj_Y     & ;\; x_j\in\Phi^-_i
\end{array}
\quad\quad\quad
g_j'\equiv\left\{
\begin{array}{ll}
\mj_{Z_q}     & ;\; x_j\in\Phi_q\neq\Phi_i \\
\mj_X         & ;\; x_j\in\Phi^+_i \\
t             & ;\; x_j\in\Phi^-_i
\end{array}\right.\right.\]
and for $1\leq j\leq n$
\[h_j\equiv\left\{
\begin{array}{ll}
\mj_{Z_q}     & ;\; y_j\in\Phi_q\neq\Phi_i \\
\mj_Y         & ;\; y_j\in\Phi^+_i \\
t             & ;\; y_j\in\Phi^-_i
\end{array}
\quad\quad\quad
h_j'\equiv\left\{
\begin{array}{ll}
\mj_{Z_q}     & ;\; y_j\in\Phi_q\neq\Phi_i \\
t             & ;\; y_j\in\Phi^+_i \\
\mj_X         & ;\; y_j\in\Phi^-_i
\end{array}\right.\right.\]
{\makebox(50,0)[bl]{}}$S(h_1,\ldots,h_n)
\alpha(Z_1,\ldots,Z_{i-1},Y,Z_{i+1},\ldots,Z_{k_\Phi})T(g_1,\ldots,g_m)=$
\\
{\makebox(115,0)[bl]{}}$S(h_1',\ldots,h_n')
\alpha(Z_1,\ldots,Z_{i-1},X,Z_{i+1},\ldots,Z_{k_\Phi})T(g_1',\ldots,g_m')$
\\[.1cm]
The equation $(\beta)$ arises when we replace $\Phi$, $m$, $n$, $T$, $\alpha$
and $S$ in $(\alpha)$ by $\Psi$, $n$, $p$, $S$, $\beta$ and $R$ respectively.
These three groups of equations are called
${\cal K}_{\Phi,\Psi}$-{\em equations}.

The following abbreviations will help us in a syntactical analysis of
the category $\cal D$. Let $[g]$ in a morphism term denote that the
morphism term $g$ may occur at that position and let $\pj_X$ denote a
composition of $q$, $q\geq 0$, morphism terms $\mj_X$. Furthermore we won't
use parentheses for composition; hence, from now on equality between
morphism terms is taken up to the associativity $(cat2)$.

\lema {\em If} $g:X\str Y$ {\em is a morphism term and} $X\in\{A,C\}$,
{\em then} $Y\in\{A,C\}$ {\em and} $g\equiv \pj_C\uf\pj_A$.
{\em In particular, if} $X\equiv C$, {\em then} $Y\equiv C$ {\em and}
$g=\mj_C$.

\dkz We proceed by induction on the complexity of the morphism term $g$.

If $g$ is a primitive morphism term, it is neither of the form
$\alpha(X_1,\ldots,X_{k_\Phi})$ nor \linebreak
$\beta(Y_1,\ldots,Y_{k_\Psi})$, since
$T(X_{\pi(x_1)},\ldots,X_{\pi(x_m)})\neq A$,
$S(Y_{\pi'(y_1)},\ldots,Y_{\pi'(y_n)})\neq C$ and $Ob({\cal D})$ is freely
generated. Hence, $g\equiv\mj_A$ or $g\equiv\mf$.

If $g$ is not primitive, then for the same reason as above, $g$ is
neither $T(g_1,\ldots,g_m)$, nor $S(h_1,\ldots,h_n)$, nor $R(t_1,\ldots,t_p)$.
Hence, $g$ is a composition $g_2g_1$ for $g_1:X\str Z$ and $g_2:Z\str Y$.
By the inductive hypothesis, since $g_1$ is of lower complexity than $g$,
$Z\in\{A,C\}$ and $g_1\equiv\pj_C\uf\pj_A$. Then by the induction hypothesis
applied to $g_2$, we have $Y\in\{A,C\}$ and $g_2\equiv\pj_C\uf\pj_A$.
Therefore $g\equiv\pj_C\uf\pj_A\pj_C\uf\pj_A$, and since $A\neq C$, we claim
$g\equiv \pj_C\uf\pj_A$. The second part of the lemma follows from the fact
that $g$ is a morphism term.
\qed
Analogously, we can prove:

\lema {\em If} $g:X\str Y$ {\em is a morphism term and} $Y\in\{A,C\}$,
{\em then} $X\in\{A,C\}$ {\em and} $g\equiv \pj_C\uf\pj_A$.

\vspace{.2cm}
Let \mtv\ abbreviate a composition of $q$, $q\geq 0$,
morphism terms of the form \linebreak
$\pj T(\pj_C\uf\pj_A,\ldots,\pj_C\uf\pj_A)\pj$, and let \msv\ and \mrv\
mean the same for $S$ and $R$ instead of $T$ respectively. Denote by
$\cal M$ the set of morphism terms of the form
\[\mrv\beta(Y_1,\ldots,Y_{k_\Psi})\msv\alpha(X_1,\ldots,X_{k_\Phi})\mtv\]
for $X_1,\ldots,X_{k_\Phi},Y_1,\ldots,Y_{k_\Psi}\in\{A,C\}$, whose type is
$T\langle A,C\rangle\str R\langle C,A\rangle$.

\lema {\em The set $\cal M$ is closed under equality.}

\dkz A substitution of equalities according to the categorial and functorial
equations doesn't change the form of a term from $\cal M$. Substitutions
of equalities according to the ``limit'' cases of $(\alpha)$ and
$(\beta)$ cause suspicion. Such is, for example, the case of substitution
according to $(\alpha)$ when $\Phi_i^+ \cap\{x_1,\ldots,x_m\}=
\Phi_i^- \cap\{y_1,\ldots,y_n\}=\emptyset$. If $g'$ is a term
obtained by such a substitution from an $\cal M$ morphism term $g$, then
an arbitrary morphism term $t:X\str C$ may occur as an argument of $T$ and
$S$, and this $X$ may occur as an argument of $\alpha$ in $g'$. However,
Lemma \arabic{section}.3 guarantees that then $X\in\{A,C\}$ and
$t\equiv\pj_C\uf\pj_A$, hence $g'$  remains in $\cal M$. We deal
with the other limit cases analogously, referring to Lemmata
\arabic{section}.2 and \arabic{section}.3 when necessary. Nonlimit cases
of substitution according to $(\alpha)$ and $(\beta)$ are obviously harmless.
\qed

\exmp Let $\Phi$ and $\Psi$ be as in Example \arabic{section}.4.
Consider the morphism term
\[h\equiv S(\mj_A^5,\mj_C^3)\alpha(A^4,C)T(\mj_A^4).\]
By g-dinaturality of $\alpha$ in the third component of $\Phi$, it is
equal to
\[h'=S(\mj^2_A,t,\mj_A^2,\mj_C^3)\alpha(A^2,X,A,C)T(\mj^3_A,t)\]
for some $t:A\str X$. Then by Lemma \arabic{section}.2, $X\in\{A,C\}$
and $t\equiv\pj_C\uf\pj_A$ which is enough for a term to remain in
$\cal M$ after the substitution of $h'$ for $h$ in it.

\vspace{.2cm}
Let $(nat)$ denote the equation
\[R\langle\mj_C,\mf\rangle\beta(C^{k_\Psi})\alpha(C^{k_\Phi})
T\langle\mf,\mj_C\rangle=
R\langle\mf,\mj_A\rangle\beta(A^{k_\Psi})\alpha(A^{k_\Phi})
T\langle\mj_A,\mf\rangle.\]
It is clear that $(nat)$ means commutativity of the diagram from
Lemma \arabic{section}.1, and therefore
\[(nat) \Leftrightarrow P(\Phi,\Psi).\]
So to prove that $P(\Phi,\Psi)$ is decidable we may use a normalization
procedure in a rewrite system corresponding to the equational theory of
$\cal M$. Actually, we have two notions of reductions. The first one is
called $CF$ ({\em categorial-functorial reduction}), and its redexes
and contracta are the following
\[\begin{array}{ccc}
{\mbox{\rm $CF$ step}} & redex & contractum\\[.1cm]
(1) & g\mj & g\\[.1cm]
(2) & \mj g & g \\[.1cm]
(3) & \quad F(h_1,\ldots,h_k)F(g_1,\ldots,g_k)\quad & F(t_1,\ldots,t_k)\\[.1cm]
(4) & F(\mj_{X_1},\ldots,\mj_{X_k}) & \mj_{F(D_1,\ldots,D_k)}
\end{array}\]
In the last two steps $F$, $k$, $g$'s, $h$'s and $t$'s satisfy the
conditions from the functorial equations above.

Since a $CF$ redex and the corresponding contractum  are equal, by Lemma
\arabic{section}.4 we have that
a term remains in $\cal M$ after a $CF$ reduction.

By the following lemma we have
that each morphism term $g$ from $\cal M$ has a unique $CF$-normal form,
which we denote by $CF(g)$.

\lema $CF$ {\em is strongly normalizing and weakly Church-Rosser.}

\dkz For strong normalization it is enough to note that a $CF$
contractum is of lower complexity than the corresponding redex.

The only interesting cases in proving that $CF$ is weakly Church-Rosser
are the following (the other cases of ramification,
roughly speaking, commute):
\begin{center}
\begin{picture}(180,160)
\put(90,20){\makebox(0,0){$F(g_{1},{\ldots},g_{k})$}}
\put(90,140){\makebox(0,0){$F(g_{1},{\ldots},g_{k})F(\mj^k)$}}
\put(20,80){\makebox(0,0){$F(g_{1}\mj,{\ldots},g_{k}\mj)$}}
\put(160,80){\makebox(0,0){$F(g_{1},{\ldots},g_{k})\mj$}}
\put(35,45){\makebox(0,0)[r]{\scriptsize$(1)$ ($k$-times) }}
\put(147,45){\makebox(0,0){\scriptsize$(1)$}}
\put(147,115){\makebox(0,0){\scriptsize$(4)$}}
\put(33,115){\makebox(0,0){\scriptsize$(3)$}}
\put(20,70){\vector(1,-1){30}}
\put(50,40){\vector(1,-1){10}}
\put(160,70){\vector(-1,-1){40}}
\put(60,130){\vector(-1,-1){40}}
\put(120,130){\vector(1,-1){40}}
\end{picture}
\end{center}
and the analogous case starting with $F(\mj^k)F(g_1,\ldots,g_k)$.

Let ${\cal M}_0$ be the set of morphism terms from $\cal M$ in
$CF$ normal form. Henceforth we use the abbreviations
$\vec X,\vec Y,\vec Z,\ldots$ for tuples of elements from the set $\{A,C\}$
and $\vec g,\vec h,\vec t,\ldots$ for tuples of elements from the set
$\{\mj_A,\mj_C,\mf\}$. From now on, a subterm in square brackets
occurs only if at least one of its arguments is $\mf$. With this notation,
we have that each member of ${\cal M}_0$ is of the shape
\[ [R(\vec t)]\beta(\vec Y)[S(\vec h)]\alpha(\vec X)[T(\vec g)]. \]

The second notion of reduction, called $D$-{\em reduction}, where $D$ stands
for dinatural, is defined on morphism terms
from ${\cal M}_0$. A peculiarity
of this reduction is that it is applicable only to the entire term as the
redex, and not to its subterms. Otherwise, it would be possible to get out of
${\cal M}_0$.

For every $i$, $1\leq i\leq k_\Phi$ and $\vec X,\vec Y, \vec g, \vec h,
\vec t$ such that both sets $\{h_j\; | \; y_j\in\Phi_i^-\}$ and
$\{g_j\; | \; x_j\in\Phi_i^+\}$ are subsets of the singleton $\{\mf\}$,
the morphism term from ${\cal M}_0$ of the following form (whose type must be
$T\langle A,C\rangle\str R\langle C,A\rangle$) 
\[ [R(\vec t)]\beta(\vec Y)[S(\vec h)]
\alpha(X_1,\ldots,X_{i-1},C,X_{i+1},\ldots,X_{k_\Phi})[T(\vec g)] \]
is a redex and
\[ [R(\vec t)]\beta(\vec Y)[S(\vec{h'})]
\alpha(X_1,\ldots,X_{i-1},A,X_{i+1},\ldots,X_{k_\Phi})[T(\vec{g'})] \]
\[{\mbox{\rm where}}\quad
g'_j\equiv\left\{
\begin{array}{ll}
g_j         & ;\; x_j\not\in\Phi_i \\
\mj_A       & ;\; x_j\in\Phi^+_i \\
\mf         & ;\; x_j\in\Phi^-_i
\end{array}
\quad{\mbox{\rm and}}\quad
h_j'\equiv\left\{
\begin{array}{ll}
h_j     & ;\; y_j\not\in\Phi_i \\
\mf     & ;\; y_j\in\Phi^+_i \\
\mj_A   & ;\; y_j\in\Phi^-_i
\end{array}\right.\right.\]
is the contractum of an $(\alpha_i)$-{\em step} of $D$ reduction.
Note that both the redex and the contractum of this step are in ${\cal M}_0$.
It follows from this fact, together with Lemmata \arabic{section}.2 and
\arabic{section}.3, that $\{h_j\; | \; y_j\in\Phi_i^+\}$ and
$\{g_j\; | \; x_j\in\Phi_i^-\}$ are subsets of $\{\mj_C\}$.

Analogously, for a fixed $1\leq i\leq k_\Psi$, we introduce 
$(\beta_i)$-{\em steps} of $D$ reduction whose redexes are terms
from ${\cal M}_0$ of the form
\[ [R(\vec t)]\beta(Y_1,\ldots,Y_{i-1},C,Y_{i+1},\ldots,Y_{k_\Psi})
[S(\vec h)]\alpha(\vec X)[T(\vec g)] \]
with both sets $\{t_j\; | \; z_j\in\Psi_i^-\}$ and
$\{h_j\; | \; y_j\in\Psi_i^+\}$ being subsets of the singleton $\{\mf\}$;
the corresponding contractum  is the morphism term
\[ [R(\vec{t'})]\beta(Y_1,\ldots,Y_{i-1},A,Y_{i+1},\ldots,Y_{k_\Psi})
[S(\vec{h'})]\alpha(\vec X)[T(\vec g)], \]
\[{\mbox{\rm where}}\quad
h'_j\equiv\left\{
\begin{array}{ll}
h_j         & ;\; y_j\not\in\Psi_i \\
\mj_A       & ;\; y_j\in\Psi^+_i \\
\mf         & ;\; x_j\in\Psi^-_i
\end{array}
\quad{\mbox{\rm and}}\quad
t_j'\equiv\left\{
\begin{array}{ll}
t_j     & ;\; z_j\not\in\Psi_i \\
\mf     & ;\; z_j\in\Psi^+_i \\
\mj_A   & ;\; z_j\in\Psi^-_i
\end{array}\right.\right.\]

\exmp For $\Phi$ and $\Psi$ as in Example \arabic{section}.4 we have the
following $(\alpha_1)$ step of $D$ reduction.
\[R(\mj_C)\beta(C^4)\alpha(C^5)T(\mf,\mj_C,\mf^2) \leadsto
R(\mj_C)\beta(C^4)S(\mf,\mj_C^7)\alpha(A,C^4)T(\mj_A,\mf^3). \]

By the following lemma we establish the uniqueness of $D$ normal form
of a morphism term from ${\cal M}_0$. We denote the $D$ normal form of $g$ by
$D(g)$.

\lema $D$ {\em is strongly normalizing and weakly Church-Rosser.}

\dkz The strong normalization property follows from the fact that every
reduction step decreases the number of $C$'s as arguments of $\alpha$
and $\beta$. For the proof that $D$ is weakly Church-Rosser, we rely on
the following facts:
\\[.2cm]
-- reduction steps $(\alpha_i)$ and $(\alpha_j)$
($(\beta_i)$ and $(\beta_j)$) commute for $i\neq j$, since connectional
components of a graph are disjoint,
\\[.1cm]
-- if a term from ${\cal M}_0$ is the redex of $(\alpha_i)$ and
$(\beta_j)$ reduction steps, then there is no $q$, $1\leq q\leq n$,
for which $y_q$ is in both $\Phi_i$ and $\Psi_j$. This is because from
the initial assumption it follows that $y_q\in\Phi_i^+$
claims $h_q\equiv\mj_C$ and $y_q\in\Psi_j^+$ claims $h_q\equiv\mf$ and
from the similar reason $y_q$ can't be a negative vertex in
$\Phi_i\cap\Psi_j$. Hence, the reduction steps $(\alpha_i)$ and $(\beta_j)$
act on disjoint sets of arguments of $T$, $S$, $R$, $\alpha$ and $\beta$
and therefore commute.
\qed

We shall find Lemmata 2.5 and 2.6 very useful for

\teo {\em Equality in} $\cal M$ {\em is decidable.}

\dkz It is enough to show that for two morphism terms $g_1$ and $g_2$ from
$\cal M$ the following equivalence holds:
\[g_1=g_2 \quad{\mbox{\rm iff}}\quad D(CF(g_1))\equiv D(CF(g_2))\]
The {\em if} part of this equivalence is trivial since all the reductions are
covered by our equations ($CF$ reductions are covered by categorial and
functorial equations and for $D$ reductions we need all ${\cal K}_{\Phi,\Psi}$
equations).

To prove the {\em only if} part, we rely on the equality axioms
(reflexivity, symmetry, transitivity and congruence), and we assume
that $g_2$ is the result of a substitution of a term for a subterm of
$g_1$ according to a ${\cal K}_{\Phi,\Psi}$ equation.
(By the equality axioms, we must have a chain of morphism terms
$g_1\equiv h_0=h_1=\ldots=h_q\equiv g_2$ such that for adjacent terms,
one is obtained from the other by a substitution described above.)
If the equation in question is a
categorial or functorial equation, then by Lemma \arabic{section}.5,
we have that $CF(g_1)\equiv CF(g_2)$; hence $D(CF(g_1))\equiv D(CF(g_2))$.
If we deal with  a dinatural equation, then it is clear that we need
just one step of $D$ reduction to reduce $CF(g_1)$ to $CF(g_2)$  or vice
versa, and therefore, by Lemma \arabic{section}.6, $D(CF(g_1))\equiv
D(CF(g_2))$.
\qed
C{\footnotesize OROLLARY}\hspace{1em} {\em The property} $P(\Phi,\Psi)$
{\em is decidable.}

Let us transform the equation $(nat)$ by deleting superfluous subterms,
if necessary, to obtain the following equation
\[(cfnat)\quad
[R\langle\mj_C,\mf\rangle]\beta(C^{k_\Psi})\alpha(C^{k_\Phi})
[T\langle\mf,\mj_C\rangle]=
[R\langle\mf,\mj_A\rangle]\beta(A^{k_\Psi})\alpha(A^{k_\Phi})
[T\langle\mj_A,\mf\rangle].\]
It is easy to see that the left-hand side ($LHS$) and the right-hand side
($RHS$) of $(cfnat)$  are in $CF$ normal form. Moreover, $RHS$ is in $D$
normal form too. Therefore, the property $P(\Phi,\Psi)$ is
equivalent to
\[D(LHS)\equiv RHS.\]
We use this equivalence in order to establish some geometrical conditions
of the amalgamation $\Phi+\Psi$, which are equivalent to $P(\Phi,\Psi)$. For
this reason we introduce the following auxiliary notation. For
a graph $\Gamma$ and $v\in V_\Gamma\setminus G_\Gamma$, let $\Gamma_v$ be
the set
$\{w\in\Gamma_{\pi(v)}\setminus G_\Gamma\; | \; \{v,w\}\not\in E_\Gamma\}$,
and let $\Gamma'_v$ be the set
$\{w\in\Gamma_{\pi(v)}\setminus G_\Gamma\; | \; \{v,w\}\in E_\Gamma\}$.
With this notation, in Example \arabic{section}.4, we have $\Psi_{y_2}=\{y_2\}$,
$\Psi'_{y_2}=\{y_3,y_7\}$, $\Phi_{x_1}=\{x_1\}$, $\Phi'_{x_1}=\{x_2,y_1\}$,
$\Phi_{x_4}=\{x_4,y_3\}$, $\Phi'_{x_4}=\emptyset$, etc.

\lema {\em For a positive} $y_i$ {\em let a morphism term from}
${\cal M}_0$ {\em in which the} $i$-{\em th argument of} $S$ {\em is}
$\mj_C$, {\em  reduce by a sequence of} $D$ {\em reductions to a
term in which this argument is} $\mf$. {\em Then this sequence of reductions
includes a step in whose
redex all the argument places from} $\Phi'_{y_i}$ {\em are occupied by}
$\mf$ {\em and the} $i$-{\em th argument of} $S$ {\em is} $\mj_C$.

\dkz Suppose that
{\scriptsize\[
[R(\vec {t^{0}})]\beta(\vec {Y^{0}})[S(\vec {h^{0}})]
\alpha(\vec {X^{0}})[T(\vec {g^{0}})]\leadsto
[R(\vec {t^{1}})]\beta(\vec {Y^{1}})[S(\vec {h^{1}})]
\alpha(\vec {X^{1}})[T(\vec {g^{1}})]\leadsto\ldots
\leadsto [R(\vec {t^{q}})]\beta(\vec {Y^{q}})[S(\vec {h^{q}})]
\alpha(\vec {X^{q}})[T(\vec {g^{q}})],\]}
\noindent is the shortest sequence of $D$ reductions for which the lemma fails. Hence,
$h_i^0\equiv\mj_C$ and $h_i^q\equiv\mf$. We claim that $h_i^1\not\equiv\mj_C$,
otherwise we would have a shorter sequence than the initial for which the
lemma fails. Also, $h_i^1$ is not $\mf$; otherwise, the first reduction step
requires all the argument places from $\Phi'_{y_i}$ in the redex to be
occupied by $\mf$, which together with $h_i^0\equiv\mj_C$ contradicts the
assumption that the lemma fails. Eventually, $h_i^1\equiv\mj_A$ is impossible
because there is no $D$ reduction step transforming $\mj_C$ to $\mj_A$
directly.
Hence, the lemma holds, since we have exhausted all the possibilities for
$h_i^1$.
\begin{flushright}$\Box$
\end{flushright}

\lema {\em For a positive} $y_i$ {\em let a morphism term from}
${\cal M}_0$ {\em in which the} $i$-{\em th argument of} $S$ {\em belongs
to the set} $\{\mf,\mj_C\}$, {\em  reduce by a sequence of} $D$
{\em reductions to a term in which this argument is} $\mj_A$.
{\em Then this sequence of reductions includes a step in whose
redex all the argument places from} $\Psi_{y_i}$ {\em are occupied by}
$\mf$ {\em and in whose contractum all the argument places from} $\Psi_{y_i}$
{\em are occupied by} $\mj_A$ {\em and all the argument places from}
$\Psi'_{y_i}$ {\em are occupied by} $\mf$.
 
\dkz Let again
{\scriptsize\[
[R(\vec {t^{0}})]\beta(\vec {Y^{0}})[S(\vec {h^{0}})]
\alpha(\vec {X^{0}})[T(\vec {g^{0}})]\leadsto
[R(\vec {t^{1}})]\beta(\vec {Y^{1}})[S(\vec {h^{1}})]
\alpha(\vec {X^{1}})[T(\vec {g^{1}})]\leadsto\ldots
\leadsto [R(\vec {t^{q}})]\beta(\vec {Y^{q}})[S(\vec {h^{q}})]
\alpha(\vec {X^{q}})[T(\vec {g^{q}})],\]}
\noindent be a shortest sequence of reductions for which the lemma fails. Note that
$h_i^0\in\{\mf,\mj_C\}$ and $h_i^q\equiv \mj_A$. Now $h_i^1$ is neither
$\mf$ nor $\mj_c$; otherwise we would have a shorter sequence for which
the lemma fails. Also, $h_i^1\not\equiv\mj_A$; otherwise, the first
reduction step requires arguments in the redex and in the contractum such that
it contradicts the assumption that the lemma fails.
\qed
We can prove the following two lemmata analogously.

\lema {\em For a negative} $y_i$ {\em let a morphism term from}
${\cal M}_0$ {\em in which the} $i$-{\em th argument of} $S$ {\em is}
$\mj_C$, {\em  reduce by a sequence of} $D$ {\em reductions to a
term in which this argument is} $\mf$. {\em Then this sequence of reductions
includes a step in whose
redex all the argument places from} $\Psi'_{y_i}$ {\em are occupied by}
$\mf$ {\em and the} $i$-{\em th argument of} $S$ {\em is} $\mj_C$.

\lema {\em For a negative} $y_i$ {\em let a morphism term from}
${\cal M}_0$ {\em in which the} $i$-{\em th argument of} $S$ {\em belongs
to the set} $\{\mf,\mj_C\}$, {\em  reduce by a sequence of} $D$
{\em reductions to a term in which this argument is} $\mj_A$.
{\em Then this sequence of reductions includes a step in whose
redex all the argument places from} $\Phi_{y_i}$ {\em are occupied by}
$\mf$ {\em and in whose contractum all the argument places from} $\Phi_{y_i}$
{\em are occupied by} $\mj_A$ {\em and all the argument places from}
$\Phi'_{y_i}$ {\em are occupied by} $\mf$.
\\[.2cm] 
In the sequel we also refer to the propositions concerning
an $x$ or a $z$ vertex
instead of $y_i$, which are analogous to the last four lemmata.

We are ready to define a geometrical criterion for $P(\Phi,\Psi)$. Let
$v_1,v_2,\ldots,v_q$ be a sequence of vertices and let
$e_1,e_2,\ldots,e_{q-1}$ be
a sequence of edges from $\Phi+\Psi$ such that $e_1=\{v_1,v_2\}$,
$e_2=\{v_2,v_3\}$, etc., and such that for each pair of adjacent
edges, one belongs to $E_\Phi$ and the other to $E_\Psi$. We call such
a pair of sequences an {\em alternating chain}. If $v_1=v_q$, then the
alternating chain is called an {\em alternating loop}.
Note that in the latter
case, the edges $e_1$ and $e_{q-1}$ are not in the same graph, and
the name alternating loop is still justified.
Here is a necessary condition for $P(\Phi,\Psi)$.

\lema {\em If} $\Phi+\Psi$ {\em provides g-dinaturality, then there are no
alternating loops in it.}

\dkz From the definition of graph it follows that the sequence of vertices
in an alternating loop consists of an even number of mutually distinct $y$'s.
Suppose now that $P(\Phi,\Psi)$ holds but that $\Phi+\Psi$ includes an
alternating loop. For the sake of clarity we use the simplest case with the
loop whose vertices are $y_i^+$ and $y_j^-$ and whose edges are
$e_1=\{y_i,y_j\}\in E_\Phi$ and $e_2=\{y_i,y_j\}\in E_\Psi$. In all the other
cases we can proceed analogously.

By the corollary of Theorem \arabic{section}.1 and by the assumption
$P(\Phi,\Psi)$ we have that the term
\[g_1\equiv[R\langle\mj_C,\mf\rangle]\beta(C^{k_\Psi})\alpha(C^{k_\Phi})
[T\langle\mf,\mj_C\rangle] \]
reduces by a D reduction to the term
\[g_2\equiv[R\langle\mf,\mj_A\rangle]\beta(A^{k_\Psi})\alpha(A^{k_\Phi})
[T\langle\mj_A,\mf\rangle]. \]
By Lemma \arabic{section}.8, this reduction must be of the form
\[g_1 \leadsto\ldots\leadsto g_3\leadsto\ldots\leadsto g_2\]
with $g_3$ an ${\cal M}_0$ morphism term whose $i$-th argument of $S$ is
$\mf$. Then by Lemma \arabic{section}.7 this reduction must be of the form
\[g_1 \leadsto\ldots\leadsto g_4\leadsto\ldots\leadsto g_3
\leadsto\ldots\leadsto g_2\]
with $g_4$ an ${\cal M}_0$ morphism term whose $j$-th argument place of $S$
is $\mf$. By Lemma \arabic{section}.9, the reduction must be of the form
\[g_1 \leadsto\ldots\leadsto g_5\leadsto\ldots\leadsto g_4
\leadsto\ldots\leadsto g_3\leadsto\ldots\leadsto g_2\]
with the $i$-th argument of $S$ being $\mf$ in $g_5$. Now we can repeat this
procedure endlessly which contradicts to the finiteness of the reduction
\qed

\noindent The necessity of our geometrical condition for $P(\Phi,\Psi)$
is of rather smaller practical interest for the purpose of proving
dinaturality of transformations. It can be used in a construction of
a countermodel for the dinaturality of composition.
However, the other direction
of the lemma above is much more useful and we are going to prove it now.
For this purposes we define the following binary relation $<_\Gamma$ in the
set of the argument places of a graph $\Gamma$: every positive left-hand
side argument place $u$ is in the relation $<_\Gamma$ with every element
of $\Gamma'_u$ and every negative right-hand side argument place $v$
is in the relation $<_\Gamma$ with every member of $\Gamma'_v$. For an
amalgamation $\Phi+\Psi$ let $<$ be the union of $<_\Phi$ and $<_\Psi$.
By this definition, we have the following chains arranged by $<$ in
Example \arabic{section}.4.
\[\begin{array}{c}
x_{1}<x_{2} \\
x_{1}<y_{1}<y_{4}<y_{5}<y_{6}<y_{8}<z_{1}\\
x_{3}<y_{2}<y_{7}<y_{8}<z_{1}\\
x_{3}<y_{2}<y_{3}\\
x_{4}
\end{array}\]

\lema {\em If there are no alternating loops in} $\Phi+\Psi$, {\em then
this amalgamation provides g-dinaturality.}

\dkz
Let $\leq$ be the reflexive and transitive closure of $<$ defined as
above in the set of argument places from $\Phi\cup\Psi$.
This set is partially ordered by $\leq$ because
of the absence of alternating loops in $\Phi+\Psi$. Suppose now that
$P(\Phi+\Psi)$ fails; i.e., for the equality $(cfnat)$ we have
\[D(LHS)\equiv [R(\vec t)]\beta(\vec Y)[S(\vec h)]\alpha(\vec X)[T(\vec g)]
\not\equiv RHS.\]
Hence, at least one of the following cases must occur in $D(LHS)$.
\begin{enumerate}
\item An argument of $R$, $S$ or $T$ is $\mj_C$.
\item For some $i$ such that $x_i$ is positive, $g_i$ is $\mf$.
\item For some $i$ such that $x_i$ is negative, $g_i$ is $\mj_A$.
\item For some $i$, $h_i$ is $\mf$.
\item For some $i$ such that $z_i$ is positive, $t_i$ is $\mj_A$.
\item For some $i$ such that $z_i$ is negative, $t_i$ is $\mf$.
\end{enumerate}

Cases 3. and 5. are impossible since the reduction preserves types of
morphism terms.

Suppose now that we have Case 1. In the ordering $\leq$, let an argument
place $v$ be minimal such that it is occupied by $\mj_C$ in $D(LHS)$.
The vertex $v$ is neither of the form $x^+$ nor $z^-$ for the same
reason as above. Suppose that
$v\equiv x_i^-$. We deal with the other cases analogously. The set
$\Phi'_{x_i}$ couldn't be empty; otherwise, $D(LHS)$ is the redex of
an $(\alpha_i)$ step of $D$ reduction. An argument place from $\Phi'_{x_i}$
couldn't be occupied by $\mj_C$ in $D(LHS)$, since for every
$v\in\Phi'_{x_i}$, $v<x_i$. If all the argument places from $\Phi'_{x_i}$
are occupied by $\mf$ in $D(LHS)$, then it is not in normal form.
If an argument place from $\Phi'_{x_i}$ is occupied by $\mj_A$ in $D(LHS)$,
then by an analogue of Lemma \arabic{section}.8 (concerning the vertex
$x_i$ instead of $y_i$) the reduction
\[LHS\leadsto\ldots\leadsto D(LHS)\]
includes a step in whose redex the $i$-th argument of $T$ is $\mf$.
Since there is no reduction transforming $\mf$ into $\mj_C$, and
since $g_i\equiv\mj_C$ in $D(LHS)$, this is impossible. Therefore, Case
1 leads to a contradiction.

Suppose now we have Case 2. As we have just seen, Case
1. doesn't obtain. If all the argument places from $\Phi_{x_i}$
are occupied by $\mf$ in $D(LHS)$, then it is not a $D$ normal form.
Let $x_j\in\Phi_{x_i}$ be occupied by $\mj_A$. The other cases are dealt with
analogously. By an analogue of Lemma \arabic{section}.8 (concerning
$x_j$ instead of $y_i$) the reduction
\[LHS\leadsto\ldots\leadsto D(LHS)\]
includes a step in whose contractum all the argument places from
$\Phi_{x_i}$ are occupied by $\mj_A$. Hence $x_i$ is occupied by $\mj_A$
in this morphism term. Since no reduction transforms $\mj_A$ into $\mf$, this
is impossible. With cases 4. and 6. we deal analogously.
\qed

Composing the previous two lemmata, we obtain the main result of the section.

\teo $P(\Phi,\Psi) \Leftrightarrow \Phi+\Psi$ {\em doesn't include alternating
loops.}

\vspace{.2cm}
Note that this theorem considers just a single component amalgamation
$\Phi+\Psi$, but as it was mentioned earlier,
this result holds universally since the notion of g-dinaturality is
defined componentwise. Also, we have reduced our considerations to functors
with arguments from one category. The generalization is trivial but
it would complicate the notation which is already, by our opinion, at the
limit of acceptability.

It is time now to compare this result with a classical one from
\cite{EK66}, which has served as an inspiration for our Theorem
\arabic{section}.2. However, the basis of \cite{EK66} (definitions of graph
and naturality) was created to fit applications involving
symmetric monoidal closed categories (cf. \cite{KM71}), and it is obvious
that we have here in mind a more involved case of
bicartesian closed categories. We
believe that our result may be applicable beyond this limitation.
It is easy to see how the part of our theorem concerning sufficiency
of the given condition for $P(\Phi,\Psi)$ captures the main result given
in \cite{EK66}.
The lack of closed curves in $\Phi+\Psi$, which was
taken there as sufficient for $P(\Phi,\Psi)$, has as a trivial consequence
the lack of alternating loops. In fact these two
conditions are equivalent in the scope of the restricted definition of graph
given in \cite{EK66}, since there are no points of ramification in $\Phi+\Psi$.
However, in our
context the presence of closed curves is harmless
for dinaturality by itself; we must instead rely on the absence of
{\em alternating} loops in amalgamations in order to guarantee dinaturality.

\section{Bicartesian closed canonical transformations}

By a {\em bicartesian closed canonical}
(also called {\em allowable}) {\em transformation} in a bicartesian
closed category $\cal B$ we mean a set of morphisms from this category
indexed by the objects from $\cal B$, defined in terms of
the special morphisms and the morphism operations from the definition given
in Section 1. Formally, it can be defined in the following manner.

Let ${\cal F}_\cb$ be the category whose objects are
functors of types $\cb^0\str\cb$, where $\cb^0$ is the trivial category
$\mj_\ast:\ast\str\ast$,
or $\cb^{l_1}\times\ldots\times\cb^{l_m}\str\cb$ for
$m\geq 0$ and $l_i\in\{-1,1\}$. We define $Ob({\cal F}_\cb)$
inductively by
\[\begin{array}{l}
1_\cb:\cb\str\cb\quad\in Ob({\cal F}_\cb),
\\[.1cm]
I:\cb^0\str\cb\; (I(\ast)=\ri)\quad\in Ob({\cal F}_\cb),
\\[.1cm]
O:\cb^0\str\cb\; (O(\ast)=\ro)\quad\in Ob({\cal F}_\cb).
\end{array}\]
If $F:\cb^{l_1}\times\ldots\times\cb^{l_m}\str\cb$ and
$G:\cb^{l_{m+1}}\times\ldots\times\cb^{l_{m+n}}\str\cb$ are in
$Ob({\cal F}_\cb)$, then
$F\otimes G:\cb^{l_1}\times\ldots\times\cb^{l_m}\times
\cb^{l_{m+1}}\times\ldots\times\cb^{l_{m+n}}\str\cb$
($F\otimes G(x_1,\ldots,x_{m+n})=F(x_1,\ldots,x_n)\times G(x_{m+1},\ldots,
x_{m+n})$),
$F\oplus G:\cb^{l_1}\times\ldots\times\cb^{l_m}\times
\cb^{l_{m+1}}\times\ldots\times\cb^{l_{m+n}}\str\cb$
($F\oplus G(x_1,\ldots,x_{m+n})=F(x_1,\ldots,x_n)+ G(x_{m+1},\ldots,
x_{m+n})$) and
$F\str G:\cb^{-l_1}\times\ldots\times\cb^{-l_m}\times
\cb^{l_{m+1}}\times\ldots\times\cb^{l_{m+n}}\str\cb$
($F\str G(x_1,\ldots,x_{m+n})=F(x_1,\ldots,x_n)\str G(x_{m+1},\ldots,
x_{m+n})$) are in $Ob({\cal F}_\cb)$.

The set of {\em canonical transformations} that we define below will be
the set of morphisms from ${\cal F}_\cb$. Each canonical transformation is a
set of $\cb$ morphisms indexed by tuples of objects from $\cb$, together
with a graph defined as in Section 2. First we define
{\em primitive canonical transformations} for every
$F:\cb^{l_1}\times\ldots\times\cb^{l_m}\str\cb$,
$G:\cb^{l_{m+1}}\times\ldots\times\cb^{l_{m+n}}\str\cb$ and
$H:\cb^{l_{m+n+1}}\times\ldots\times\cb^{l_{m+n+p}}\str\cb$ from
$Ob({\cal F}_\cb)$.

$\mj_F=\{\mj_{F(\vec{A})}\; |\; \vec{A}\in (Ob(\cb))^m\}$
is a primitive canonical transformation from $F$ to $F$
whose graph consists of
vertices $x_1^{l_1},\ldots,x_m^{l_m},y_1^{l_1},\ldots, y_m^{l_m}$ and
edges $\{x_1,y_1\}$, $\ldots$, $\{x_m,y_m\}$.

$\md_F=\{\md_{F(\vec{A})}\; |\; \vec{A}\in (Ob(\cb))^m\}$ is a
primitive canonical transformation from $F\otimes I$ to $F$ whose graph is
identical to the graph of $\mj_F$.

$\mc_{F,G}=\{\mc_{F(\vec{A}),G(\vec{B})}\; |\; \vec{A}\in (Ob(\cb))^m,
\vec{B}\in (Ob(\cb))^n\}$ is a primitive canonical transformation from
$F\otimes G$ to $G\otimes F$ with the graph whose vertices are
$x_1^{l_1},\ldots,x_m^{l_m}$, $x_{m+1}^{l_{m+1}}$, $\ldots,$
$x_{m+n}^{l_{m+n}}$,
$y_1^{l_{m+1}},$ $\ldots,$ $y_n^{l_{m+n}},$
$y_{n+1}^{l_1},$ $\ldots,$ $y_{n+m}^{l_m}$ and
whose edges are $\{x_1,y_{n+1}\},\ldots,\{x_m,y_{n+m}\}$, $\{x_{m+1},y_1\}$,
$\ldots,$  $\{x_{m+n},y_n\}$.

$\mw_F=\{\mw_{F(\vec{A})}\; |\; \vec{A}\in (Ob(\cb))^m\}$ is a primitive
canonical transformation from $F$ to $F\otimes F$ whose graph consists of
vertices
$x_1^{l_1},\ldots,x_m^{l_m}$, $y_1^{l_{1}},\ldots, y_m^{l_{m}}$,
$y_{m+1}^{l_1},\ldots,y_{m+m}^{l_m}$ and edges $\{x_1,y_1\}$,
$\{x_1,y_{m+1}\}$, $\ldots,$ $\{x_m,y_m\}$, $\{x_m,y_{2m}\}$.

$\mk_F=\{\mk_{F(\vec{A})}\; |\; \vec{A}\in (Ob(\cb))^m\}$ is a primitive
canonical transformation from $F$ to $I$ with the graph whose vertices
are $x_1^{l_1},\ldots,x_m^{l_m},g_1,\ldots,g_m$ and whose edges are
$\{x_1,g_1\},\ldots,\{x_m,g_m\}$.

$\mep_{F,G}=
\{\mep_{F(\vec{A}),G(\vec{B})}\; |\; \vec{A}\in (Ob(\cb))^m,\vec{B}\in
(Ob(\cb))^n\}$ is a primitive canonical transformation from
$F\otimes(F\str G)$ to $G$ with the graph whose vertices are
$x_1^{l_1},\ldots,x_m^{l_m}$, $x_{m+1}^{-l_{1}},\ldots,x_{m+m}^{-l_{m}}$,
$x_{2m+1}^{l_{m+1}},\ldots,x_{2m+n}^{l_{m+n}}$,
$y_1^{l_{m+1}},\ldots, y_n^{l_{m+n}}$ and whose edges are $\{x_1,x_{m+1}\}$,
$\ldots$, $\{x_m,x_{2m}\}$, $\{x_{2m+1},y_1\}$, $\ldots$, $\{x_{2m+n},y_n\}$.

Analogously, we define the primitive canonical transformations
$\md^i_F$ from $F$ to $F\otimes \ri$,
$\mb^\str_{F,G,H}$ from $F\otimes(G\otimes H)$ to
$(F\otimes G)\otimes H$,
$\mb^\rts_{F,G,H}$ from $(F\otimes G)\otimes H$ to
$F\otimes (G\otimes H)$, $\ml_F$ from $\ro$ to $F$,
$\ml^1_{F,G}$ from $F$ to $F\oplus G$, $\ml^2_{F,G}$ from $G$ to $F\oplus G$,
$\mm_F$ from $F\oplus F$ to $F$ and $\met_{F,G}$ from $G$ to
$F\str(F\otimes G)$ with corresponding graphs. It is not difficult to show
that every primitive canonical transformation is g-dinatural with respect
to the associated graph.

Next we define the following operations on canonical transformations.

If $\alpha=\{\alpha(A_1,\ldots,A_{k_\Gamma})\;|\;A_1,\ldots,A_{k_\Gamma}
\in Ob(\cb)\}$ is a canonical transformation from $F$ to $G$ with
the graph $\Gamma$, then for $l\geq 1$
\[\alpha^{i_1,\ldots,i_l}=\{\alpha(A_1,\ldots,A_{k_\Gamma})\;
|\;A_1,\ldots,A_{k_\Gamma}
\in Ob(\cb), A_{i_1}=A_{i_2}=\ldots=A_{i_l}\}\]
is a canonical transformation from $F$ to $G$ with the graph obtained from
$\Gamma$ by addition of edges between the vertices from the components
$i_1,\ldots,i_l$ in order to obtain one component of the new graph. We call
$\alpha^{i_1,\ldots,i_l}$ a {\em subtransformation} of $\alpha$. It is easy
to verify that if $\alpha$ and $\beta$ are canonical transformation from
$F$ to $G$ and if $\beta\subset\alpha$, then $\beta$ is a subtransformation
of $\alpha$. Also, if a canonical transformation is g-dinatural, then each of
its subtransformations is g-dinatural, too.

If 
$\alpha=\{\alpha(A_1,\ldots,A_{k_\Phi})\;|\;A_1,\ldots,A_{k_\Phi}\in Ob(\cb)\}$
and
$\beta=\{\beta(B_1,\ldots,B_{k_\Psi})\;|\;B_1,\ldots,B_{k_\Psi}\in Ob(\cb)\}$
are two canonical transformations from $F$ to $G$ and from $H$ to $J$
respectively, then
\[\begin{array}{l}
\alpha\otimes\beta=\{\alpha(A_1,\ldots,A_{k_\Phi})\times
\beta(B_1,\ldots,B_{k_\Psi})\;|\;A_1,\ldots,B_{k_\Psi}\in Ob(\cb)\},
\\[.1cm]
\alpha\oplus\beta=\{\alpha(A_1,\ldots,A_{k_\Phi})+
\beta(B_1,\ldots,B_{k_\Psi})\;|\;A_1,\ldots,B_{k_\Psi}\in Ob(\cb)\},
\\[.1cm]
\alpha\str\beta=\{\alpha(A_1,\ldots,A_{k_\Phi})\str
\beta(B_1,\ldots,B_{k_\Psi})\;|\;A_1,\ldots,B_{k_\Psi}\in Ob(\cb)\}
\end{array}\]
are canonical transformations from $F\otimes H$ to $G\otimes J$,
from $F\oplus H$ to $G\oplus J$, and from $G\str H$ to $F\str J$
respectively.
If $\Phi$ is the graph of $\alpha$ and $\Psi$ is the graph of $\beta$,
then the graphs of $\alpha\otimes\beta$, $\alpha\oplus\beta$ and
$\alpha\str\beta$ are obtained as disjoint unions of $\Phi$ and $\Psi$, where in
the last case, $\Phi$ occurs inverted. We denote these graphs by
$\Phi\otimes\Psi$, $\Phi\oplus\Psi$ and $\Phi\str\Psi$ respectively.

\exmp Let $\Phi$ be the graph on left-hand side and $\Psi$ the graph on
right-hand side of the picture below.

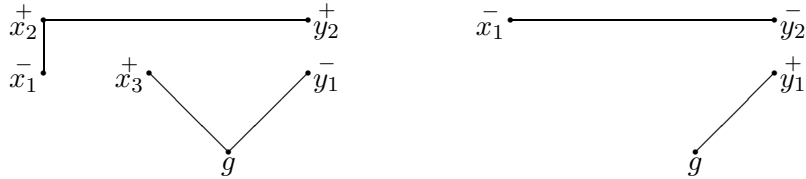
\begin{center}
\begin{picture}(140,70)(0,10)
\put(20,50){\circle*{2}}
\put(18,50){\makebox(0,0)[r]{$\stackrel{-}{x_{1}}$}}
\put(20,70){\circle*{2}}
\put(18,70){\makebox(0,0)[r]{$\stackrel{+}{x_{2}}$}}
\put(60,50){\circle*{2}}
\put(58,50){\makebox(0,0)[r]{$\stackrel{+}{x_{3}}$}}
\put(90,20){\circle*{2}}
\put(90,18){\makebox(0,0)[t]{$g$}}
\put(120,50){\circle*{2}}
\put(122,50){\makebox(0,0)[l]{$\stackrel{-}{y_{1}}$}}
\put(120,70){\circle*{2}}
\put(122,70){\makebox(0,0)[l]{$\stackrel{+}{y_{2}}$}}
\put(20,50){\line(0,1){20}}
\put(20,70){\line(1,0){100}}
\put(60,50){\line(1,-1){30}}
\put(90,20){\line(1,1){30}}
\end{picture}
\quad\quad\quad
\begin{picture}(140,70)(0,10)
\put(20,70){\circle*{2}}
\put(18,70){\makebox(0,0)[r]{$\stackrel{-}{x_{1}}$}}
\put(90,20){\circle*{2}}
\put(90,18){\makebox(0,0)[t]{$g$}}
\put(120,50){\circle*{2}}
\put(122,50){\makebox(0,0)[l]{$\stackrel{+}{y_{1}}$}}
\put(120,70){\circle*{2}}
\put(122,70){\makebox(0,0)[l]{$\stackrel{-}{y_{2}}$}}
\put(20,70){\line(1,0){100}}
\put(90,20){\line(1,1){30}}
\end{picture}
\end{center}
Then $\Phi\otimes\Psi$ and $\Phi\oplus\Psi$ are identical and given by the
diagram on left-hand side and $\Phi\str\Psi$ is given by the diagram on
right-hand side below.
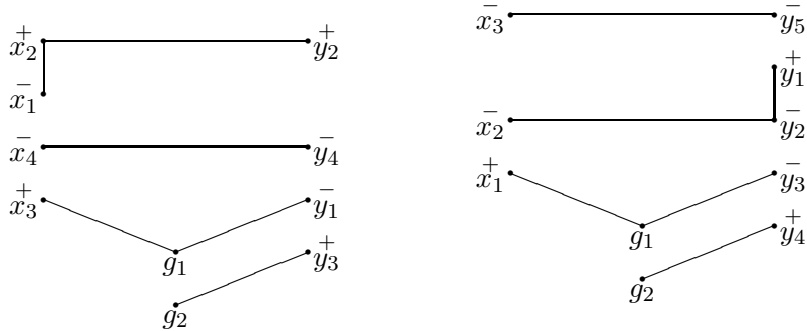
\begin{center}
\begin{picture}(140,120)(0,10)
\put(20,60){\circle*{2}}
\put(18,60){\makebox(0,0)[r]{$\stackrel{+}{x_{3}}$}}
\put(20,80){\circle*{2}}
\put(18,80){\makebox(0,0)[r]{$\stackrel{-}{x_{4}}$}}
\put(20,100){\circle*{2}}
\put(18,100){\makebox(0,0)[r]{$\stackrel{-}{x_{1}}$}}
\put(20,120){\circle*{2}}
\put(18,120){\makebox(0,0)[r]{$\stackrel{+}{x_{2}}$}}
\put(70,20){\circle*{2}}
\put(70,18){\makebox(0,0)[t]{$g_{2}$}}
\put(70,40){\circle*{2}}
\put(70,38){\makebox(0,0)[t]{$g_{1}$}}
\put(120,40){\circle*{2}}
\put(122,40){\makebox(0,0)[l]{$\stackrel{+}{y_{3}}$}}
\put(120,60){\circle*{2}}
\put(122,60){\makebox(0,0)[l]{$\stackrel{-}{y_{1}}$}}
\put(120,80){\circle*{2}}
\put(122,80){\makebox(0,0)[l]{$\stackrel{-}{y_{4}}$}}
\put(120,120){\circle*{2}}
\put(122,120){\makebox(0,0)[l]{$\stackrel{+}{y_{2}}$}}
\put(20,100){\line(0,1){20}}
\put(20,120){\line(1,0){100}}
\put(20,80){\line(1,0){100}}
\put(20,60){\line(5,-2){50}}
\put(70,40){\line(5,2){50}}
\put(70,20){\line(5,2){50}}
\end{picture}
\quad\quad\quad
\begin{picture}(140,130)
\put(20,60){\circle*{2}}
\put(18,60){\makebox(0,0)[r]{$\stackrel{+}{x_{1}}$}}
\put(20,80){\circle*{2}}
\put(18,80){\makebox(0,0)[r]{$\stackrel{-}{x_{2}}$}}
\put(20,120){\circle*{2}}
\put(18,120){\makebox(0,0)[r]{$\stackrel{-}{x_{3}}$}}
\put(70,20){\circle*{2}}
\put(70,18){\makebox(0,0)[t]{$g_{2}$}}
\put(70,40){\circle*{2}}
\put(70,38){\makebox(0,0)[t]{$g_{1}$}}
\put(120,40){\circle*{2}}
\put(122,40){\makebox(0,0)[l]{$\stackrel{+}{y_{4}}$}}
\put(120,60){\circle*{2}}
\put(122,60){\makebox(0,0)[l]{$\stackrel{-}{y_{3}}$}}
\put(120,80){\circle*{2}}
\put(122,80){\makebox(0,0)[l]{$\stackrel{-}{y_{2}}$}}
\put(120,100){\circle*{2}}
\put(122,100){\makebox(0,0)[l]{$\stackrel{+}{y_{1}}$}}
\put(120,120){\circle*{2}}
\put(122,120){\makebox(0,0)[l]{$\stackrel{-}{y_{5}}$}}
\put(120,80){\line(0,1){20}}
\put(20,120){\line(1,0){100}}
\put(20,80){\line(1,0){100}}
\put(20,60){\line(5,-2){50}}
\put(70,40){\line(5,2){50}}
\put(70,20){\line(5,2){50}}
\end{picture}
\end{center}

Moreover, canonical transformations are closed
under composition defined as in Section 2; i.e.
if $\alpha$ is a canonical transformation from $F$ to $G$ with the graph
$\Phi$ and $\beta$ is a canonical transformation from $G$ to $H$ with the graph
$\Psi$, then $\beta\alpha$ defined as in Section 2, is a canonical
transformation with the graph $\Psi\Phi$.

It is easy to verify that ${\cal F}_\cb$ is a category with the identity morphism
for $F$ being $\mj_F$ and the composition of $\alpha$ and $\beta$ being $\beta\alpha$
defined as above. We leave the details about the structure of this category for
another occasion.

Our aim is to show that all the morphisms from ${\cal F}_\cb$ are
g-dinatural transformations. It is easy to see that the only obstacle for this is the
composition of canonical transformations. To show that composition is
now harmless too,
we use the results from Section 2 and the procedure of cut elimination in an
adequate sequent system.

The following example shows that the results from Section 2 are not sufficient
for our aims before a further analysis of properties peculiar to
bicartesian closed categories.

\exmp Let $\alpha$ be the canonical transformation obtained from the following
composition of canonical transformations (from now on, we associate compositions
to the right)
\[(\mep_{1,1}\otimes\mj_{1\str 1})\mb^\str_{1,(1\str 1),(1\str 1)}
(\mep_{1,1}\otimes\mj_{(1\str 1)\otimes(1\str 1)})
\mb^\str_{1,1\str 1, (1\str1)\otimes(1\str 1)}
(\mj_1\otimes(\mj_{1\str 1}\otimes\mw_{1\str 1}))(\mj_1\otimes\mw_{1\str 1})\]
and let $\beta$ be $\mep_{1,1}$.
From the facts that the primitive canonical transformations are g-dinatural, that
$\otimes$ preserves g-dinaturality, and from Theorem 2.2, it follows that $\alpha$
and $\beta$ are g-dinatural transformations whose amalgamation of graphs is given
by the following diagram.
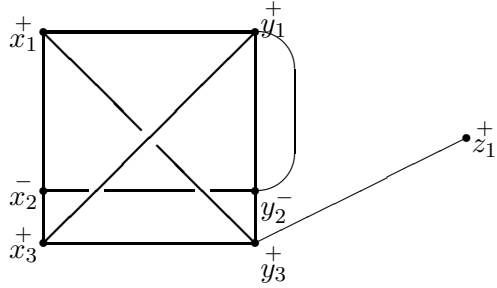
\begin{center}
\begin{picture}(180,100)
\put(10,10){\circle*{3}}
\put(8,10){\makebox(0,0)[r]{$\stackrel{+}{x_{3}}$}}
\put(90,10){\circle*{3}}
\put(92,10){\makebox(0,0)[tl]{$\stackrel{+}{y_{3}}$}}
\put(10,30){\circle*{3}}
\put(8,30){\makebox(0,0)[r]{$\stackrel{-}{x_{2}}$}}
\put(90,30){\circle*{3}}
\put(92,30){\makebox(0,0)[tl]{$y^-_{2}$}}
\put(10,90){\circle*{3}}
\put(8,90){\makebox(0,0)[r]{$\stackrel{+}{x_{1}}$}}
\put(90,90){\circle*{3}}
\put(92,90){\makebox(0,0)[bl]{$\stackrel{+}{y_{1}}$}}
\put(170,50){\circle*{3}}
\put(172,50){\makebox(0,0)[l]{$\stackrel{+}{z_{1}}$}}
{\thicklines
\put(10,10){\line(1,0){80}}
\put(10,10){\line(0,1){80}}
\put(90,90){\line(-1,0){80}}
\put(90,90){\line(0,-1){80}}
\put(90,90){\line(-1,-1){80}}
\put(10,90){\line(1,-1){37}}
\put(90,10){\line(-1,1){37}}
\put(10,30){\line(1,0){17}}
\put(33,30){\line(1,0){34}}
\put(90,30){\line(-1,0){17}}}
\put(90,10){\line(2,1){80}}
\put(90,60){\oval(30,60)[r]}
\end{picture}
\end{center}
Since an alternating loop occurs in this amalgamation, by Theorem 2.2.
there is a composition of g-dinatural transformations with such graphs,
which is not g-dinatural. Of course, it doesn't mean that $\beta\alpha$ is
not g-dinatural. However, each element of $\beta\alpha$ is in the composition
of canonical transformations
\[\alpha_1=\mb^\str_{1\otimes(1\str 1),(1\str 1),(1\str 1)}
\mb^\str_{1,1\str 1, (1\str1)\otimes(1\str 1)}
(\mj_1\otimes(\mj_{1\str 1}\otimes\mw_{1\str 1}))(\mj_1\otimes\mw_{1\str 1})\]
and
\[\beta_1=\mep_{1,1}(\mep_{1,1}\otimes\mj_{1\str 1})
(\mep_{1,1}\otimes\mj_{(1\str 1)\otimes(1\str 1)})\]
which in turn gives that $\beta\alpha$ is a subtransformation of $\beta_1\alpha_1$.
The g-dinaturality of $\beta_1$, and of $\beta_1\alpha_1$ too, follows
from Theorem 2.2. Hence, from these two facts it follows that $\beta\alpha$ is
g-dinatural. In the sequel, we generalize the idea from the example above
to the
case of an arbitrary composition of canonical transformations. For this
purpose
we need the following definitions.

Let $\alpha$ be a canonical transformation. Denote by $C(\alpha)$ the set
of canonical transformations defined inductively by
\\[.2cm]
-- $\alpha\in C(\alpha)$,
\\[.1cm]
-- if $\beta\in C(\alpha)$ and $F\in Ob({\cal F}_\cb)$, then
$(\beta\otimes\mj_F)$ and $(\mj_F\otimes \beta)$ are in $C(\alpha)$.
\\[.2cm]
Let $\xi_{F,G,H}$ from $F\otimes(G\oplus H)$ to
$(F\otimes G)\oplus(F\otimes H)$
be the following canonical transformation.
\[\mep_{F,(F\otimes G)\oplus(F\otimes H)}(\mj_F\otimes
\mm_{F\str((F\otimes G)\oplus(F\otimes H))}
((\mj_F\str\ml^1_{F\otimes G,F\otimes H})\met_{F,G}+
((\mj_F\str\ml^2_{F\otimes G,F\otimes H})\met_{F,H})))\]

Next we define the set \ccb\ of {\em constructible} canonical transformations.
This name comes from the analogous notion from \cite{KM71}.
\begin{enumerate}
\item Primitive canonical transformations are in \ccb.
\item If $\alpha$ from $F$ to $G$ is in $C(\beta)$ for $\beta$ be
among $\mb^\rts_{T,S,R}$,
$\mb^\str_{T,S,R}$, $\mc_{T,S}$, $\mw_T$, $\mk_T$, $\md_T$, $\md^i_T$ for
some $T,S,R\in Ob({\cal F}_\cb)$, and $\gamma$ from $G$ to $H$ is in \ccb,
then $\gamma\alpha$ is in \ccb.
\item If $\alpha$ and $\beta$ are in \ccb\ then $\alpha\otimes\beta$ is in
\ccb.
\item If $\alpha$ from $F\otimes G$ to $H$ and $\beta$ from $J\otimes G$
to $H$ are in \ccb, then \linebreak $\mm_H(\alpha\oplus\beta)(\mc_{G,F}\oplus
\mc_{G,J})\xi_{F,J,G}\mc_{F\oplus J,G}$ is in \ccb.
\item If $\alpha$ from $F$ to $G$ is in \ccb\ then $\ml^1_{G,H}\alpha$ and
$\ml^2_{G,H}\alpha$ are in \ccb.
\item If $\alpha$ from $F\otimes G$ to $H$ is in \ccb\ then
$(\mj_F\str\alpha)\met_{F,G}$ is in \ccb.
\item If $\alpha$ from $F$ to $G$ and $\beta$ from $H\otimes J$ to $T$
are in \ccb\, then $\beta((\mep_{G,H}(\alpha\otimes\mj_{G\str H}))\otimes
\mj_J)$ is in \ccb.
\end{enumerate}

\lema {\em Each constructible canonical transformation is g-dinatural with
respect to its own graph.}

\dkz It is easy to verify that the primitive canonical transformations are
g-dinatural. (This follows from the equations ($\md$), ($\mb$), ($\mc$),
($\mw$), ($\mk$), ($\ml^1$), ($\ml^2$), ($\ml$), ($\mm$), ($\mep 1$),
($\met 1$), ($\mep 2$) and ($\met 2$).) Also it is easy to see that if
$\alpha$ and $\beta$ are g-dinatural, then such are $\alpha\otimes \beta$,
$\alpha\oplus\beta$ and $\alpha\str\beta$, too. For the rest, we rely on
Theorem 2.2.

\section{A Category-like Sequent System for Intuitionistic
Propositional Logic}

In this section we carry out a cut elimination procedure in an
auxiliary sequent system for intuitionistic propositional logic, which will
help us in dealing with the dinaturality of bicartesian closed canonical
transformations.

This sequent system, which we call, \J\ is introduced as follows.
Let \F\ be generated from a countable set $\cal L$, whose members we
call {\em propositional letters}, with the constants $\top$ and $\bot$ and
the binary connectives $\I$, $\IL$ and $\str$. We call the members of \F\
{\em formulae}, and use the schematic letters $A,B,C,\ldots,A_1,\ldots$ for
them.
{\em Sequents} of \J\ are of the form $A \vdash B$ for $A$ and $B$
in \F. We call $A$ in $A \vdash B$ the {\em antecedent}, and $B$ the
{\em consequent} of the sequent. In order to introduce the rules of
inference of \J\ we need
the following auxiliary notion of {\em $\I$-context}, which corresponds to
the notion of (poly)functor in categories. A $\I$-context is defined
inductively as follows:
\begin{description}
\item[$1^\circ$] the symbol $\Box$ is a $\I$-context,
\item[$2^\circ$] if $G$ is a $\I$-context and $A \in \cal F$,
then $(G \I A)$ and $(A \I G)$ are $\I$-contexts,
\item[$3^\circ$] if $G$ and $H$ are $\I$-contexts, then $(G \I H)$ is a
$\I$-context.
\end{description}
For a $\I$-context $F$ we say that it is a $\I_{1}$-context if the symbol
$\Box$ occurs in $F$ exactly once. For $G$ a $\I$-context and $A \in \cal F$,
we obtain $G(A)$ by substituting $A$ for $\Box$ in $G$, e.g., if
$F \equiv (B \I \Box) \I C$, then $F(A)=(B \I A) \I C$.
\\[0.3cm]
The {\em axioms} of \J\ are
\[ a_A : A \vdash A,\quad\Pi_A:\bot\vdash A, {\mbox{\hspace{2ex}for every
\hspace{1ex}}}A \in \cal F, \]
\\[0.3cm]
The {\em structural rules} of \J\ are
\[
\begin{array}{ll}
(\beta_{F}^{\leftarrow})\rza
\f{F (A \I (B \I C)) \vdash D}
{F ((A \I B) \I C)
\vdash D} \rzb
&
(\beta_{F}^{\str})\rza
\f{F ((A \I B) \I C) \vdash D}
{F(A \I (B \I C)) \vdash D}
\\
(\gamma_F)\rza \f{F(A \I B) \vdash C}
{F(B \I A) \vdash C}
&
\\
(\omega_F)\rza
\f{F(A \I A) \vdash B}{F(A) \vdash B}
&
(\theta^{A}_{F})\rza
\f{F(\top) \vdash B}{F(A) \vdash B}
\\
(\tau_F)\rza
\f{F(A) \vdash B}{F(A \I \top) \vdash B}
&
(\tau^i_F)\rza
\f{F(A \I \top) \vdash B}{F(A) \vdash B}
\\
(\mix{G})\rza
\f{A \vdash B \rzb G(B) \vdash C}{G(A) \vdash C},
\end{array}
\]
where $F$ is a $\I_1$ context and $G$ is a $\I$ context. The last rule is
called {\em mix} and we refer to it by ($\mix{}$) when the context $G$ is
irrelevant.
\\[0.3cm]
The {\em rules for connectives} are 
\[
\begin{array}{ll}
(\I)\rza
\f{A \vdash C \rzb B \vdash D}{A \I B \vdash C \I D}\rzb
&
\\
(\Diamond)\rza
\f{A\I C\vdash D\rzb B\I C\vdash D}{(A\IL B)\I C\vdash D}\rzb
&
(\d{C})\rza
\f{A\vdash B}{A\vdash B\IL C}\rzb
(\l{B})\rza
\f{A\vdash C}{A\vdash B\IL C}
\\
(\ast)\rza
\f{A \I B \vdash C}{B \vdash A \str C}
&
(\rhd)\rza
\f{A \vdash B \rzb C \I D \vdash E}{(A \I (B \str C)) \I D \vdash E}
\end{array}
\]

A {\em proof} of a sequent $A \vdash B$ in \J\ is a binary tree
with sequents in its nodes, such that $A \vdash B$ is in the root,
axioms are in the leaves and consecutive nodes are connected by
some of the inference rules above.

It is not difficult to see that the underlying
logic of \J\ is intuitionistic propositional logic.
The differences between \J\ and Gentzen's system $LJ$ introduced in
\cite{G35} are that in \J\ we have just one
meta-logical symbol ($\vdash$) in the sequents: we omit Gentzen's commas
in the antecedents, whose role is now covered by the logical connective $\I$.
We can't have empty either the antecedent or the consequent of a
sequent in \L. The logical constant $\top$ serves to fill gaps in antecedents.
These discrepancies between \J\ and $LJ$ arise because in \J\ we want
antecedents and consequents of sequents to be of the same sort (namely
members of $\cal F$) and this enables us to look at an \J\ sequent as
an arrow with the source being the antecedent and the target the consequent
of the sequent.

Our ($\I$) is a rule of simultaneous introduction of
the connective $\I$ on the both sides of a sequent: there is no
a counterpart for this rule in $LJ$. This difference is not categorially
motivated though it emphasizes functoriality of the connective $\I$.
We also believe that \J\ completely separates
structural rules from the rules for connectives. On the other hand,
the $LJ$ rules \&-IS and \&-IA (see 1.22. of \cite{G35}) have hidden
interchanges, contractions and thinnings.

Since we prove the cut-elimination theorem through elimination of {\em mix},
as Gentzen did too, we have postulated the mix rule ($\mix{}$) as primitive.
However, this mix is something different from Gentzen's mix. It is
liberal in the sense that the $\I$-context $G$ in ($\mix{G}$) need not
to capture all factors $B$ (see the definition below) as arguments in
$G(B)$. This means that the formula $B$ may be a factor of $A$ in Step $2^\circ$ of the
construction of the $\I$-context $G$; i.e. mix need not to ``swallow''
all the occurrences of $B$ in $G(B)$. There are no categorial
reasons to prefer cut to such a mix. In both cases,
we don't have categorial composition of arrows corresponding to both
premises of the rule, but a more involved composition of the right premise
with an image of the left premise under the functor corresponding
to a $\I$-context. The only difference is that in the case of cut this is
always a $\I_1$-context.

An advantage of \J\ is that its proofs can be easily coded.
For example the proof
\[\f{\afrac{p \vdash p} \rzb \f{q \vdash q}{q \I \top \vdash q}}
{\f{(p \I (p \str q)) \I \top \vdash q}{p \I (p \str q) \vdash q}}\]
is coded by
\[ \tau^i_{\Box} (a_p \rhd \tau_{\Box} a_q) \]
This fact helps when we want to postulate equalities that should hold
between the proofs of \J.

For $G$ a $\I$-context and $\pi$ a proof,
we denote by $G(\pi)$ the proof coded by the term obtained from $G$
after the substitution $a_A$ for every $A$ and the code of $\pi$ for
every $\Box$ in $G$.

For the proof of the main result of this section we need the following
notions of
{\em degree} and {\em rank}. The degree of a formula is the number of logical
connectives in it. However, because of the categorially motivated elimination
of the comma, the symbol $\I$ plays a double
role and in order to define rank, we
define as follows a set of {\em factors} of $A$, for every $A \in \cal F$:
\begin{description}
\item[$1^\circ$] $A$ is a factor of $A$,
\item[$2^\circ$] if $A$ is of the form $A_1 \I A_2$ then every factor of
$A_1$ or $A_2$ is a factor of $A$.
\end{description}

Now, we introduce (in the style of Do\v sen) an auxiliary
indexing of consequents
and factors of antecedents in a mixless proof of \J\, which will help us
in defining the rank of an occurrence of a formula in such a proof.
First we index all the consequents and all the factors of antecedents of
axioms by 1 and inductively proceed as follows. In all the structural
rules and the rule ($\rhd$) the index of the consequent in the conclusion
is increased by 1. In ($\Diamond$) the index of the consequent in the
conclusion is the maximum of the two indices of consequents of both premises
increased by 1.
In ($\I$), (\d{C}), (\l{C}) and ($\ast$) the index of the consequent in the
conclusion is 1. Every factor of the antecedent preserved by a rule
has the index increased by 1, and all the factors introduced by the rule
have index 1 in the conclusion. In $(\omega_F)$ the
occurrence of $A$ in the conclusion is indexed by the maximum of the
indices of distinguished $A$'s in the premise, increased by 1.
In the example of
the proof given above this indexing looks like
\[\f{\afrac{p^1 \vdash p^1} \rzb \f{q^1 \vdash q^1}
{(q^2 \I \top^1)^1 \vdash q^2}}
{\f{((p^2 \I (p \str q)^1)^1 \I \top^2)^1 \vdash q^3}
{(p^3 \I (p \str q)^2)^2 \vdash q^4}}\]
Then the rank of an occurrence of a formula in a mixless proof is given
by its index.

The following theorem corresponds to Gentzen's {\em Hauptsatz} of \cite{G35}.

\teo {\em  Every proof in \J\ can be transformed into a proof of the same
root-sequent with no applications of the rule} ($\mix{}$).
\\[0.2cm]
\dkz
As in the standard cut-elimination procedure it is enough to consider
a proof $\pi$ whose last rule is ($\mix{G}$) for a $\I$-context $G$,
and there is no more application
of ($\mix{}$) in $\pi$. So let our proof be of the form
\[\f{\fp{\pi_1}{A \vdash B} \rzb \fp{\pi_2}{G(B) \vdash C}}{G(A) \vdash C}
\pravilo{\mix{G}}\]
with $\pi_1$ and $\pi_2$ mixless. Then we define the {\em degree of this proof}
as the degree of $B$ and the {\em rank of this proof} as the sum of the
{\em left rank}, i.e. the rank
of the occurrence of $B$ in the left premise of $\mix{G}$,
in the subproof $\pi_1$, and the {\em right rank}, i.e. the maximum of
all ranks of distinguished
factors $B$ in the right premise of $\mix{G}$ in the subproof $\pi_2$.
Then we prove our theorem by induction on the lexicographically
ordered pairs $\langle d, r \rangle$ for the degree $d$ and the rank $r$
of the proof.
\\[0.3cm]
1. $r=2$
\\[0.2cm]
1.1. $\pi_1$ or $\pi_2$ are axioms
\\[.2cm]
1.1.1. Suppose $\pi$ is of the form
\[\f{\afrac{a_B:B\vdash B} \rzb \fp{\pi_2}{G(B)\vdash C}}{G(B)\vdash C}\]
Then this proof is transformed into the proof
\[\fp{\pi_2}{G(B)\vdash C}\]
which is mixless.
\\[.2cm]
1.1.2. If $\pi$ is of the form
\[\f{\fp{\pi_1}{A\vdash B} \rzb \afrac{a_{G(B)}:G(B)\vdash G(B)}}
{G(A)\vdash G(B)}\]
Then this proof is transformed into the proof
\[\fp{G(\pi_1)}{G(A)\vdash G(B)}\]
which is of course mixless.
\\[.2cm]
1.1.3. If $\pi$ is of the form
\[\f{\afrac{\Pi_B:\bot\vdash B} \rzb \fp{\pi_2}{G(B)\vdash C}}
{G(\bot)\vdash C}\]
Then this proof is transformed into the proof of the form
\[\Bpak{\Pi_C:\bot\vdash C}{$(\tau),(\gamma),(\theta)$}{G(\bot)\vdash C}\]
\\[.2cm]
1.1.4. Finally, if $\pi$ is of the form
\[\f{\fp{\pi_1}{A\vdash \bot} \rzb \afrac{\Pi_C:\bot \vdash C}}{A\vdash C}\]
Then, since the left rank of this proof is 1, $A$ must be $\bot$ and
$\pi$ is transformed into $\Pi_C:\bot \vdash C$.
\\[.2cm]
1.2. $\pi_1$ ends with ($\I$)

Suppose $\pi$ is of the form
\[\f{\f{\fp{\pi'_1}{A_1 \vdash B_1} \rzb \fp{\pi''_1}{A_2 \vdash B_{2}}}
{A_1 \I A_2 \vdash B_1 \I B_2}\pravilo{$\I$} \rzc \fp{\pi_2}
{G(B_1 \I B_2) \vdash C}}
{G(A_1 \I A_2) \vdash C}\pravilo{$\mix{G}$} \]
Then this proof is transformed into the proof
\[\f{\fp{\pi'_1}{A_1 \vdash B_1} \rzc
\f{\fp{\pi''_1}{A_2 \vdash B_2} \rzb \fp{\pi_2}{G(B_1\I B_2) \vdash C}}
{G(B_1\I A_2) \vdash C}\pravilo{$\mix{}$}}
{G(A_1 \I A_2) \vdash C}\pravilo{$\mix{}$} \]
where both applications of ($\mix{}$) have lower degree.
\\[.2cm]
1.3. $\pi_1$ ends with ($\ast$)
\\[.2cm]
1.3.1. $\pi_2$ ends with ($\theta$)

Suppose $\pi$ is of the form
\[\f{\f{\fp{\pi'_1}{B_1\I A\vdash B_2}}{A\vdash B_1\str B_2}\pravilo{$\ast$}
\rzb\quad\quad \f{\fp{\pi'_2}{F(\top)\vdash C}}{G(B_1\str B_2)\vdash C}
\pravilo{$\theta^{G_1(B_1\str B_2)}_F$}}{G(A)\vdash C}\]
Then this proof is transformed into the proof
\[\f{\fp{\pi'_2}{F(\top)\vdash C}}{G(A)\vdash C}
\pravilo{$\theta^{G_1(A)}_F$}\]
\vspace{.2cm}1.3.2. $\pi_2$ ends with ($\rhd$)

Suppose $\pi$ is of the form
\[\f{\f{\fp{\pi'_1}{B_1\I A\vdash B_2}}{A\vdash B_1\str B_2}\pravilo{$\ast$}
\rzb\quad\quad \f{\fp{\pi'_2}{D\vdash B_1} \rzb \fp{\pi''_2}{B_2\I E\vdash C}}
{(D\I (B_1\str B_2))\I E\vdash C}
\pravilo{$\rhd$}}{(D\I A)\I E\vdash C}\]
Then this proof is transformed into the proof
\[\f{\fp{\pi'_2}{D\vdash B_1} \rzb \f{\fp{\pi'_1}{B_1\I A\vdash B_2} \rzb
\fp{\pi''_2}{B_2\I E\vdash C}}
{(B_1\I A)\I E\vdash C}\pravilo{$\mix{}$}}{(D\I A)\I E\vdash C}
\pravilo{$\mix{}$}\]
with both applications of ($\mix{}$) of the lower degree.
\\[.2cm]
1.4. $\pi_1$ ends with (\d{$B_2$}), or analogously with (\l{$B_1$})
\\[.2cm]
1.4.1. $\pi_2$ ends with ($\theta$) is analogous to 1.3.1.
\\[.2cm]
1.4.2. $\pi_2$ ends with ($\Diamond$)

Suppose $\pi$ is of the form
\[\f{\f{\fp{\pi'_1}{A\vdash B_1}}{A\vdash B_1\IL B_2}\pravilo{\d{$B_2$}}
\rzb\quad\quad \f{\fp{\pi'_2}{B_1\I D\vdash C} \rzb \fp{\pi''_2}
{B_2\I D\vdash C}}{(B_1\IL B_2)\I D\vdash C}
\pravilo{$\Diamond$}}{A\I D\vdash C}\]
Then this proof is transformed into the proof 
\[\f{\fp{\pi'_1}{A\vdash B_1} \rzb \fp{\pi'_2}{B_1\I D\vdash C}}{A\I D\vdash C}
\pravilo{$\mix{}$}\]
with the smaller degree.
\\[.2cm]
2. $r>2$
\\[.2cm]
2.1. the right rank is $>1$
\\[.2cm]
2.1.1. $\pi_2$ ends with a structural rule ($\sigma$), i.e.,
$\pi$ is of the form
\[\f{\fp{\pi_1}{A\vdash B} \rzb \f{\fp{\pi_2}{G_1(B)\vdash C}}
{G(B)\vdash C}\pravilo{$\sigma$}}{G(A)\vdash C}\]
\vspace{.2cm}2.1.1.1. If all the distinguished $B$'s in the right premise of $\mix{G}$ in
$\pi$ have indices grater than 1 (by 2.1, at least one such  $B$ must occur)
then this proof is transformed into the proof
\[\f{\fp{\pi_1}{A\vdash B} \rzb \fp{\pi_2}{G_1(B)\vdash C}}
{\f{G_1(A)\vdash C}{G(A)\vdash C}\pravilo{$\sigma$}}
\pravilo{\mix{G_1}}\]
whose subproof ending with \mix{G_1} has the rank lower by 1.
\\[.2cm]
2.1.1.2. If one of the distinguished $B$'s in the right premise of $\mix{G}$
in $\pi$ is indexed by 1 (note that except for ($\theta$), in the conclusion
of a structural rule, every formula has at most one occurrence indexed by 1),
then $\pi$ is transformed into the proof
\[\f{\afrac{\fp{\pi_1}{A\vdash B}} \rzb \f{\fp{\pi_1}{A\vdash B} \rzb
\fp{\pi_2}{G_1(B)\vdash C}}
{\f{G_1(A)\vdash C}{F(B)\vdash C}\pravilo{$\sigma$}}\pravilo{\mix{G_1}}}
{G(A)\vdash C}
\pravilo{\mix{F}}\]
for a $\I_1$-context $F$ (except when $\sigma$ is an application of
($\theta$) in which case $F$ is a $\I$-context) such that $F(A)\equiv G(A)$.
In this proof, the subproof ending with the upper mix has the rank decreased
by 1, and the right rank of the lower mix remains equal to 1 after
the elimination of the upper mix.
\\[.2cm]
2.1.2. $\pi_2$ ends with ($\I$)

Suppose $\pi$ is of the form
\[\f{\fp{\pi_1}{A\vdash B} \rzb \f{\fp{\pi'_2}{G_1(B)\vdash C_1} \rzb
\fp{\pi''_2}{G_2(B)\vdash C_2}}{G(B)\vdash C_1\I C_2}\pravilo{$\I$}}
{G(A)\vdash C_1\I C_2}\]
Then this proof is transformed into the proof
\[\f{\f{\fp{\pi_1}{A\vdash B} \rzb \fp{\pi'_2}{G_1(B)\vdash C_1}}{G_1(A)\vdash C_1}
\pravilo{\mix{G_1}} \rzb\quad\quad
\f{\fp{\pi_1}{A\vdash B} \rzb \fp{\pi''_2}{G_2(B)\vdash C_2}}
{G_2(A)\vdash C_2}\pravilo{\mix{G_2}}}{G(A)\vdash C_1\I C_2}\pravilo{$\I$}\]
in which both subproofs ending with \mix{G_1} and \mix{G_2} are of the
lower ranks. There is also a simplified variant of 2.1.2 with no
distinguished $B$'s in the antecedent of a premise of the rule ($\I$).

In all the cases below, the subproofs of the reduced proofs ending
with the applications of (\mix{}), have a smaller rank than $\pi$.
\\[.2cm]
2.1.3. $\pi_2$ ends with ($\ast$)

Suppose $\pi$ is of the form
\[\f{\fp{\pi_1}{A\vdash B} \rzb \f{\fp{\pi'_2}{C_1\I G(B)\vdash C_2}}
{G(B)\vdash C_1\str C_2}\pravilo{$\ast$}}
{G(A)\vdash C_1\str C_2}\]
Then this proof is transformed into the proof
\[\f{\fp{\pi_1}{A\vdash B} \rzb \fp{\pi_2}{C_1\I G(B)\vdash C_2}}
{\f{C_1\I G(A)\vdash C_2}{G(A)\vdash C_1\str C_2}\pravilo{$\ast$}}
\pravilo{\mix{}}\]
\vspace{.2cm}2.1.4. $\pi_2$ ends with ($\rhd$)
\\[.2cm]
2.1.4.1. Suppose that $\pi$ is of the form
\[\f{\fp{\pi_1}{A\vdash B} \rzb\quad\quad \f{\fp{\pi'_2}{G_1(B)\vdash B_1}
\rzb \fp{\pi''_2}{B_2\I G_2(B)\vdash C}}{(G_1(B)\I B)\I G_2(B)\vdash C}
\pravilo{$\rhd$}}{(G_1(A)\I A)\I G_2(A)\vdash C}\]
Than this proof is transformed into the proof
\[\f{\fp{\pi_1}{A\vdash B} \rzb \f{\f{\fp{\pi_1}{A\vdash B} \rzb
\fp{\pi'_2}{G_1(B)\vdash B_1}}{G_1(A)\vdash B_1}\pravilo{\mix{}}
\rzb\quad\quad \f{\fp{\pi_1}{A\vdash B} \rzb
\fp{\pi''_2}{B_2\I G_2(B)\vdash C}}{B_2\I G_2(A)\vdash C}\pravilo{\mix{}}}
{(G_1(A)\I B)\I G_2(A)\vdash C}\pravilo{$\rhd$}}
{(G_1(A)\I A)\I G_2(A)\vdash C}\pravilo{\mix{}}\]
\vspace{.2cm}2.1.4.2. Suppose that $\pi$ is of the form
\[\f{\fp{\pi_1}{A\vdash B} \rzb\quad\quad \f{\fp{\pi'_2}{B_1\vdash B_2}
\rzb \fp{\pi''_2}{B_3\I G_1(B)\vdash C}}{B\I G_1(B)\vdash C}
\pravilo{$\rhd$}}{A\I G_1(A)\vdash C}\]
Than this proof is transformed into the proof
\[\f{\fp{\pi_1}{A\vdash B} \rzb \f{\fp{\pi'_2}{B_1\vdash B_2} \rzb
\f{\fp{\pi_1}{A\vdash B} \rzb \fp{\pi''_2}{B_3\I G_1(B)\vdash C}}
{B_3\I G_1(A)\vdash C}\pravilo{\mix{}}}{B\I G_1(A)\vdash C}\pravilo{$\rhd$}}
{A\I G_1(A)\vdash C}\pravilo{\mix{}}\]
\vspace{.2cm}2.1.4.3. Suppose that $\pi$ is of the form
\[\f{\fp{\pi_1}{A\vdash B} \rzb\quad\quad \f{\fp{\pi'_2}{G_1(B)\vdash D}
\rzb \fp{\pi''_2}{E\I G_2(B)\vdash C}}{(G_1(B)\I (D\str E))\I G_2(B)\vdash C}
\pravilo{$\rhd$}}{(G_1(A)\I (D\str E))\I G_2(A)\vdash C}\]
Than this proof is transformed into the proof
\[\f{\f{\fp{\pi_1}{A\vdash B} \rzb \fp{\pi'_2}{G_1(B)\vdash D}}
{G_1(A)\vdash D}\pravilo{\mix{}} \rzb\quad\quad
\f{\fp{\pi_1}{A\vdash B} \rzb \fp{\pi''_2}{E\I G_2(B)}}{E\I G_2(A)\vdash C}
\pravilo{\mix{}}}{(G_1(A)\I (D\str E))\I G_2(A)\vdash C}\pravilo{$\rhd$}\]
There are also simplified variants of 2.1.4.1. and 2.1.4.3. with no
distinguished $B$'s in $G_1$ or $G_2$ which we won't discuss here separately.
\\[.2cm]
2.1.5. $\pi_2$ ends with ($\Diamond$)
\\[.2cm]
2.1.5.1. Suppose that $\pi$ is of the form
\[\f{\fp{\pi_1}{A\vdash B} \rzb\quad\quad \f{\fp{\pi'_2}{B_1\I G_1(B)\vdash C}
\rzb \fp{\pi''_2}{B_2\I G_1(B)\vdash C}}{B\I G_1(B)\vdash C}
\pravilo{$\Diamond$}}{A\I G_1(A)\vdash C}\]
Then this proof is transformed into the proof
\[\f{\fp{\pi_1}{A\vdash B} \rzb \f{\f{\fp{\pi_1}{A\vdash B}\rzb
\fp{\pi'_2}{B_1\I G_1(B)\vdash C}}{B_1\I G_1(A)\vdash C}\pravilo{\mix{}}
\rzb\quad\quad \f{\fp{\pi_1}{A\vdash B} \rzb \fp{\pi''_2}
{B_2\I G_1(B)\vdash C}}{B_2\I G_1(A)\vdash C}\pravilo{\mix{}}}
{B\I G_1(A)\vdash C}\pravilo{$\Diamond$}}{A\I G_1(A)\vdash C}
\pravilo{\mix{}}\]
\vspace{.2cm}2.1.5.2. Suppose that $\pi$ is of the form
\[\f{\fp{\pi_1}{A\vdash B} \rzb\quad\quad \f{\fp{\pi'_2}{D_1\I G_1(B)\vdash C}
\rzb \fp{\pi''_2}{D_2\I G_1(B)\vdash C}}{(D_1\IL D_2)\I G_1(B)\vdash C}
\pravilo{$\Diamond$}}{(D_1\IL D_2)\I G_1(A)\vdash C}\]
Then this proof is transformed into the proof
\[\f{\f{\fp{\pi_1}{A\vdash B} \rzb \fp{\pi'_2}{D_1\I G_1(B)\vdash C}}
{D_1\I G_1(A)\vdash C}\pravilo{\mix{}} \rzb\quad\quad
\f{\fp{\pi_1}{A\vdash B} \rzb \fp{\pi''_2}{B_2\I G_1(B)\vdash C}}
{B_2\I G_1(A)\vdash C}\pravilo{\mix{}}}{(D_1\IL D_2)\I G_1(A)\vdash C}
\pravilo{$\Diamond$}\]
\vspace{.2cm}2.1.6. $\pi_2$ ends with (\d{$C_2$})

Suppose that $\pi$ is of the form
\[\f
{\fp {\pi_1} {A\vdash B} \rzb \f {\fp {\pi_2} {G(B)\vdash C_1}}
{G(B)\vdash C_1\IL C_2}\pravilo{\d{$C_2$}}}{G(A)\vdash C_1\IL C_2}\]
Then this proof is transformed into the proof
\[\f
{\fp {\pi_1} {A\vdash B} \rzb \fp {\pi_2} {G(B)\vdash C_1}}
{\f {G(A)\vdash C_1} {G(A)\vdash C_1\IL C_2} \pravilo{\d{$C_2$}}}
\pravilo{\mix{}}\]
The case of (\l{$C_1$}) instead of (\d{$C_2$}) is dealt with analogously.
\\[.2cm]
2.2. The right rank is 1 and the left rank is greater than 1.
\\[.2cm]
If $\pi_2$ is the axiom $a_{G(B)}$, then we proceed as in 1.1.2. If
$\pi_2$ ends with an application of ($\theta$), then we proceed as in
1.3.1. In all the remaining cases $G$ must be a $\I_1$-context
\\[.2cm]
2.2.1. $\pi_1$ ends with a structural rule

Suppose that $\pi$ is of the form
\[\f{\f{\fp{\pi'_1}{A_1\vdash B}}{A\vdash B}\pravilo{$\sigma$}
\rzb\quad\quad \fp{\pi_2}{G(B)\vdash C}}{G(A)\vdash C}\]
Then this proof is transformed into the proof
\[\f{\fp{\pi'_1}{A_1\vdash B} \rzb \fp{\pi_2}{G(B)\vdash C}}
{\f{G(A_1)\vdash C}{G(A)\vdash C}\pravilo{$\sigma$}}\pravilo{\mix{}}\]
\vspace{.2cm}2.2.2. $\pi_1$ ends with ($\rhd$)

Suppose that $\pi$ is of the form
\[\f{\f{\fp{\pi'_1}{A_1\vdash A_2} \rzb \fp{\pi''_1}{A_3\I A_4\vdash B}}
{(A_1\I (A_2\str A_3))\I A_4\vdash B}\pravilo{$\rhd$} \rzb\quad\quad
\fp{\pi_2}{G(B)\vdash C}}{G((A_1\I (A_2\str A_3))\I A_4)\vdash C}\]
Then this proof is transformed into the proof
\[\f
{\bfrac{\fp{\pi'_1}{A_1\vdash A_2}}{18pt} \rzb
\f{\fp{\pi''_1}{A_3\I A_4\vdash B} \rzb \fp{\pi_2}{G(B)\vdash C}}
{\Bpak{G(A_3\I A_4)\vdash C}{($\beta$), ($\gamma$)}{A_3\I D\vdash C}}
\pravilo{\mix{}}}
{\Bpak{(A_1\I (A_2\str A_3))\I D\str C}{($\beta$), ($\gamma$)}
{G((A_1\I (A_2\str A_3))\I A_4)\vdash C}}\pravilo{$\rhd$}\]
\vspace{.2cm}2.2.3. Eventually, if $\pi_1$ ends with ($\Diamond$) and
$\pi$ is of the form
\[\f{\f{\fp{\pi'_1}{A_1\I A_3\vdash B}\rzb \fp{\pi''_1}{A_2\I A_3\vdash B}}
{(A_1\IL A_2)\I A_3\vdash B}\pravilo{$\Diamond$} \rzb\quad\quad
\fp{\pi_2}{G(B)\vdash C}}{G((A_1\IL A_2)\I A_3)\vdash C}\]
Then this proof is transformed into the proof
\[\f{\f{\fp{\pi'_1}{A_1\I A_3\vdash B} \rzb \fp{\pi_2}{G(B)\vdash C}}
{\Bpak{G(A_1\I A_3)\vdash C}{($\beta$), ($\gamma$)}{A_1\I D\vdash C}}\pravilo{\mix{}}
\rzb\quad\quad \f{\fp{\pi''_1}{A_2\I A_3\vdash B} \rzb \fp{\pi_2}{G(B)\vdash C}}
{\Bpak{G(A_2\I A_3)\vdash C}{($\beta$), ($\gamma$)}{A_2\I D\vdash C}}
\pravilo{\mix{}}}
{\Bpak{(A_1\IL A_2)\I D\vdash C}{($\beta$), ($\gamma$)}
{G((A_1\IL A_2)\I D)\vdash C}}\pravilo{$\Diamond$}\]
\qed

\section{The embedding of $\cal J$ into a free bicartesian closed category}

Let \bcc\ be the bicartesian closed category freely generated by
the set of objects $\cal L$ used in Section 4. The morphisms of this category
can be viewed as equivalence classes of {\em morphism terms} generated from
$\mj_A$, $\md_A$, $\md^i_A$,
$\mb^\str_{A,B,C}$, $\mb^\rts_{A,B,C}$, $\mc_{A,B}$, $\mw_A$, $\mk_A$,
$\mep_{A,B}$, $\met_{A,B}$, $\ml_A$, $\ml^1_{A,B}$, $\ml^2_{A,B}$ and
$\mm_A$ for some objects $A,B,C$ of \bcc\ with the operations $\times$,
$+$, $\str$ and $\circ$, modulo bicartesian closed equations given in
Section 1.

Now we define translations from the set of $\cal J$-formulae and the set of
$\cal J$-proofs to $Ob(\bcc)$ and the set of morphism terms, respectively.
Denote both these translations by $t$.

Let $t$ be the identity on $\cal L$ and inductively defined as follows.
(In the following definition, $\mfv$ is a naturally extracted functor from the
$\I$-context $F$, and the indices of special morphisms can be easily
reconstructed.)
\[\begin{array}{lll}
t(\top)=\ri, & t(\bot)=\ro & \\[.1cm]
t(A\I B)=t(A)\times t(B), & t(A\IL B)=t(A)+ t(B), & t(A\str B)=t(A)\str t(B),
\end{array}\]
\[\begin{array}{ll}
t(a_A)=\mj_{t(A)}, & t(\Pi_A)=\ml_{t(A)},
\\[.1cm]
t(\beta^\str_F(\pi))=t(\pi)\circ\mfv(\mb^\str), &
t(\beta^\rts_F(\pi))=t(\pi)\circ\mfv(\mb^\rts),
\\[.1cm]
t(\gamma_F(\pi))=t(\pi)\circ\mfv(\mc), &
\\[.1cm]
t(\omega_F(\pi))=t(\pi)\circ\mfv(\mw), &
t(\theta^A_F(\pi))=t(\pi)\circ\mfv(\mk_A),
\\[.1cm]
t(\tau_F(\pi))=t(\pi)\circ\mfv(\md), & t(\tau^i_F(\pi))=t(\pi)\circ\mfv(\md^i),
\\[.1cm]
t(\pi_2\:\mix{G}\:\pi_1)=t(\pi_2)\circ\mgv(t(\pi_1)), &
\\[.1cm]
t(\pi_1\I\pi_2)=t(\pi_1)\times t(\pi_2), &
\\[.1cm]
t(\pi_1\Diamond\pi_2)=\mm\circ(t(\pi_1)+t(\pi_2))\circ(\mc+\mc)\circ\xi\circ
\mc, &
\\[.1cm]
t(\pi^{\d{C}})=\ml^1_{\_\:,C}\circ t(\pi), & t(\pi^{\l{C}})=\ml^2_{C,\:\_}\circ t(\pi),
\\[.1cm]
t(\pi^\ast)=(\mj\str t(\pi))\circ \met, &
\\[.1cm]
t(\pi_1\rhd\pi_2)=t(\pi_2)\circ((\mep\circ(t(\pi_1)\times\mj))\times\mj).
\end{array}\]

The translation $t'$ that is inverse to $t$ on the set $Ob(\bcc)$ is
defined on the set of morphism terms as follows. (Here we write $A'$ instead of
$t'(A)$ and $f'$ instead of $t'(f)$.)
\[\begin{array}{ll}
t'(\mj_A)=a_{A'}, & t'(\ml_A)=\Pi_{A'},
\\[.1cm]
t'(\mb^\str_{A,B,C})=\beta^\str_\Box a_{A'\I(B'\I C')}, &
t'(\mb^\rts_{A,B,C})=\beta^\rts_\Box a_{(A'\I B')\I C')},
\\[.1cm]
t'(\mc_{A,B})=\gamma_\Box a_{B'\I A'}, &
\\[.1cm]
t'(\mw_A)=\omega_\Box a_{A'\I A'}, &
t'(\mk_A)=\theta^{A'}_\Box a_\top,
\\[.1cm]
t'(\md_A)=\tau_\Box a_{A'}, & t'(\md^i_A)=\tau^i_\Box a_{A'\I \top},
\\[.1cm]
t'(\ml^1_{A,B})=(a_{A'})^{\d{B'}}, & t'(\ml^2_{A,B})=(a_{B'})^{\l{A'}},
\\[.1cm]
t'(\mm_A)=\tau^i_\Box(\tau_\Box a_{A'}\Diamond\tau_\Box a_{a'}), &
\\[.1cm]
t'(\mep_{A,B})=\tau^i_\Box(a_{A'}\rhd\tau_\Box a_{B'}), &
t'(\met_{A,B})=(a_{A'\I B'})^\ast,
\\[.1cm]
t'(f\times g)=f'\I g', & t'(f\str g)=(\tau^i_\Box(f'\rhd\tau_\Box g'))^\ast,
\\[.1cm]
t'(f+g)=\tau^i_\Box((\tau_\Box (f')^{\d{D'}})\Diamond(\tau_\Box(g')^{\l{B'}})), &
\\[.1cm]
t'(g\circ f)=g'\: \mix{\Box}\: f'. &
\end{array}\]

\lema {\em For every morphism term} $g$, $t(t'(g))=g$.

\dkz By induction on the complexity of $g$.

\lema {\em In each step of our cut elimination procedure, which transforms}
$\pi$ {\em into} $\pi'$,
{\em we have} $t(\pi)=t(\pi')$ in \bcc.

\dkz Long, tedious but more or less trivial. In steps where we
were not precise about the order of application of rules in the
transformed proof we rely on some coherece properties, like for example
in Case 1.1.3, we use the
fact that $t(G(\bot))$ is isomorphic to \ro\ and therefore the
order of application of $(\tau)$, $(\gamma)$ and $(\theta)$ is
arbitrary.

\lema {\em For every morphism term}  $g$, {\em there is a mixless
proof} $\pi$ {\em of} $\cal J$, {\em such that} $g=t(\pi)$.

\dkz Let $\pi_1$ be $t'(g)$ and let $\pi$ be the mixless proof obtained
from $\pi_1$ by our cut elimination procedure. Then by Lemma
\arabic{section}.15,
$t(\pi_1)=t(\pi)$, and by Lemma \arabic{section}.14, $g=t(t'(g))=t(\pi_1)=t(\pi)$.
\qed

Since the mixless proofs of $\cal J$ correspond to the constructible canonical
transformations, we can derive the following lemma.

\lema {\em Every canonical transformation
from an arbitrary bicartesian closed category is a subtransformation of a
constructible canonical transformation.}

\dkz By Lemma \arabic{section}.16 and the universal property of
\bcc\ it follows that
for every canonical transformation $\alpha$ from $F$ to $G$ there
exists a constructible canonical transformation $\beta$ from $F$ to $G$
such that each member of $\alpha$ is equal to a member of $\beta$.
By our definition, this fact is sufficient for $\alpha$ being a
subtransformation of $\beta$.
\qed

From Lemmata 3.13 and 5.17 we have the following.

\teo {\em Every bicartesian closed canonical transformation is g-dinatural.}
\\[.2cm]
From this theorem and the remark after Example 2.3, it follows
that every bicartesian closed canonical transformation is dinatural
in the classical sense. Moreover, one has to bear in mind that this property
is provable regardless of the choice of language for bicartesian closed
categories.

Our proof covers a result from \cite{GSS92} where the authors have used a
normalization in a natural deduction system for the fragment of
intuitionistic propositional logic that corresponds to cartesian
closed categories, to show that all canonical transformations from
these categories are dinatural. Since there are still some difficulties
with normalization in clumsy $\lambda$-calculuses for full
intuitionistic propositional logic, we find an advantage in
sequent systems, which are sufficient to deal with the
questions of dinaturality.
The definitions of operations on objects in the
underlying functor category given in \cite{GSS92} are different from our
operations in
${\cal F}_\cb$, and a consequence of this difference is that the functor
category of \cite{GSS92} is cartesian closed, whereas our ${\cal F}_\cb$ is
just symmetric monoidal closed.

Investigations of dinaturality are often tied to investigations of
coherence. Some results (cf. \cite{B93}) claim that this connection is
very strict. However, our graphs, though appropriate for dinaturality, are
inadequate for coherence. We leave all this questions about coherence
for another occasion.
\\[.5cm]
{\bf Acknowledgments.} Most of these results are from the author's
Ph.D. thesis, written under the direction of Professor Kosta Do\v sen,
to whom the author is grateful very much.

\end{document}